\documentclass[10pt]{article}
\usepackage{graphicx}
\usepackage{amssymb}
\newtheorem{theorem}{Theorem}
\newtheorem{proposition}[theorem]{Proposition}
\newtheorem{lemma}[theorem]{Lemma}
\newtheorem{corollary}[theorem]{Corollary}

\input tcilatex

\begin{document}

\begin{titlepage}

\vskip0.5truecm

\vskip1.0truecm

\begin{center}

{\LARGE \bf Area-preserving diffeomorphisms of the torus whose rotation sets have non-empty interior}

\end{center}

\vskip  0.4truecm

\centerline {{\large Salvador Addas-Zanata}}

\vskip 0.2truecm

\centerline { {\sl Instituto de Matem\'atica e Estat\'\i stica }}
\centerline {{\sl Universidade de S\~ao Paulo}}
\centerline {{\sl Rua do Mat\~ao 1010, Cidade Universit\'aria,}} 
\centerline {{\sl 05508-090 S\~ao Paulo, SP, Brazil}}
 
\vskip 0.7truecm

\begin{abstract}

In this paper we consider $C^{1+\epsilon}$ area-preserving diffeomorphisms of the 
torus 
$f,$ either homotopic to the identity or to Dehn twists.
We suppose that $f$ has a lift $\widetilde{f}$ to the plane such that its 
rotation set has interior and prove, among other 
things that if zero is an interior point of the rotation set, then
there exists a hyperbolic $\widetilde{f}$-periodic point 
$\widetilde{Q}$$\in {\rm I}\negthinspace {\rm R^2}$ such that
$W^u(\widetilde{Q})$ intersects $W^s(\widetilde{Q}+(a,b))$ for 
all integers $(a,b)$, which implies that 
$\overline{W^u(\widetilde{Q})}$ is invariant under integer translations. Moreover, 
$\overline{W^u(\widetilde{Q})}=\overline{W^s(\widetilde{Q})}$ and
$\widetilde{f}$ restricted to $\overline{W^u(\widetilde{Q})}$ is invariant and
 topologically mixing.   
Each connected component of the complement of $\overline{W^u(\widetilde{Q})}$ is
a disk with diameter uniformly bounded from above. If $f$ is transitive, then 
$\overline{W^u(\widetilde{Q})}=$${\rm I}\negthinspace {\rm R^2}$ and 
$\widetilde{f}$ is topologically mixing in the whole plane.

\end{abstract} 

\vskip 0.3truecm

\vskip 2.0truecm


\vskip 0.8truecm

\noindent{\bf e-mail:} sazanata@ime.usp.br

\vskip 1.0truecm


\vfill
\hrule
\noindent{\footnotesize{The author is partially supported 
by CNPq, grant: 304803/06-5}}

\end{titlepage}

\baselineskip=6.2mm

\section{Introduction and main results}

One of the most well understood chapters of dynamics of surface
homeomorphisms is the case of the torus. Any orientation preserving
homeomorphism $f$ of the torus can be associated in a canonical way with a
two by two matrix $A$ with integer coefficients and determinant one.
Depending on this matrix, there are basically three types of maps\ (if we
are allowed to consider iterates of $f$ and coordinate changes):

\begin{enumerate}
\item  for some integer $n>0,$ $A^n$ is the identity and in this case, $f^n$
is said to be homotopic to the identity;

\item  for some integer $n>0,$ $A^n$ is, up to a conjugation, equal to $%
\left( 
\begin{array}{cc}
1 & k \\ 
0 & 1
\end{array}
\right) ,$ for some integer $k\neq 0.$ If this is the case, $f^n$ is said to
be homotopic to a Dehn twist;

\item  $A$ is hyperbolic, that is, $A$ has real eigenvalues $\lambda $ and $%
\mu $ and $\left| \mu \right| <1<\left| \lambda \right| ;$
\end{enumerate}

This paper concerns to the first two cases, namely we suppose throughout the
paper that the matrix $A$ associated to $f$ is either the identity or it is
equal to $\left( 
\begin{array}{cc}
1 & k \\ 
0 & 1 
\end{array}
\right) ,$ for some integer $k\neq 0.$ In these cases, it is possible to
consider a rotation set which, roughly speaking, ''measures'' how orbits in
the torus rotate with respect to the homology (below we will present precise
definitions). It is a generalization of the rotation number of an
orientation preserving circle homeomorphism to this two-dimensional context
and a lot of work has been done on this subject. For instance, one wants to
show connections between dynamical properties of $f$ and geometric
properties of the rotation set, see \cite{franksrat}, \cite{misiu} and \cite
{llibre} for maybe the first references on this problem and also which sets
can be realized as rotation sets of a torus homeomorphism, see for instance 
\cite{kwaerg} and \cite{kwanon}.

Our initial motivation was the following question: Is the ''complexity'' of
the map $f,$ in some sense, shared by its lift to the plane? The objective
was to look at an area-preserving homeomorphism $f$ of the torus with a
rotation set with non-empty interior. And then give some consequences of
these hypotheses to the lift of $f$ to the plane, denoted $\widetilde{f}:%
{\rm I\negthinspace R^2\rightarrow I\negthinspace R^2,}$ used to compute the
rotation set. In some cases, we could assume other hypothesis on $f,$ for
instance like being transitive and we wanted to know what happened with $
\widetilde{f}.$ Unfortunately, to prove something interesting, we had to
leave the continuous world and assume $C^{1+\epsilon }$ differentiability
(for some $\epsilon >0$), so that we could use Pesin theory. Apart from our
main theorems, this paper has some lemmas that may have interest by
themselves. For instance, it is known since \cite{franksrat}, \cite{doeff}
and \cite{eu4} that rational points in the interior of the rotation set are
realized by periodic orbits. One of our results says that given a rational
rotation vector in the interior of the rotation set, there is a hyperbolic
periodic point which realizes this rotation vector and it has, what Kwapisz
and Swanson \cite{kwabar} called a rotary horseshoe, that is, the union of
its stable and unstable manifolds contains a homotopically non-trivial
simple closed curve in the torus. This sort of property is useful to prove
most of the results in this paper. In particular, given two rational vectors
in the interior of the rotation set, there are hyperbolic periodic points
which realize these rotation vectors and the stable manifold of one of them
intersects the unstable manifold of the other periodic point (and vice
versa) in a topologically transverse way (see definition 9). 

We have similar results for maps homotopic to the identity and homotopic to
Dehn twists. Although some proofs are slightly different, most of them work
in both cases, with just simple adjustments. Only when the ideas involved
are different, we will present separate proofs in the Dehn twist case.
Studying maps homotopic to Dehn twists is certainly of great interest, since
for instance the well-known Chirikov standard map, $S_{M,k}:{\rm %
T^2\rightarrow T^2,}$ given by 
\begin{equation}
\label{standardmap}S_{M,k}(x,y)=(x+y+k\sin (2\pi x)\text{ }mod1,y+k\sin
(2\pi x)\text{ }mod1), 
\end{equation}
where $k>0$ is a parameter, is homotopic to $(x,y)\rightarrow $ $\left( 
\begin{array}{cc}
1 & 1 \\ 
0 & 1 
\end{array}
\right) \left( 
\begin{array}{c}
x \\ 
y 
\end{array}
\right) (mod1)^2.$ And, maybe some of the hardest questions in surface
dynamics refer to the standard map. For instance, one of them, which has
some relation with the results proved in this paper is the following:

\begin{description}
\item[Conjecture]  : Is there $k>0$ such that $S_{M,k}$ is transitive?
\end{description}

This is an very difficult question. Just to give an idea, Pedro Duarte \cite
{duarte} proved that for a residual set of large values of $k,$ $S_{M,k}$
has lots of elliptic islands, therefore it is not transitive.

%
%
%
%
%
%
%
%
%
%

Finally, we shall point out that in \cite{faand} and \cite{faand1} results
with more or less the same flavor as the ones from this paper have been
shown to be true for area preserving torus  homeomorphisms  homotopic to the
identity, using different techniques.

Before presenting our results, we need some definitions.

\vskip0.2truecm

{\large {\bf Definitions:}}

\begin{enumerate}
\item  Let ${\rm T^2}={\rm I}\negthinspace {\rm R^2}/{\rm Z\negthinspace
\negthinspace Z^2}$ be the flat torus and let $p:{\rm I}\negthinspace {\rm %
R^2}\longrightarrow {\rm T^2}$ and $\pi :{\rm I}\negthinspace {\rm R^2}%
\longrightarrow S^1\times {\rm I}\negthinspace {\rm R}$ be the associated
covering maps. Coordinates are denoted as $(\widetilde{x},\widetilde{y})\in 
{\rm I}\negthinspace {\rm R^2,}$ $(\widehat{x},\widehat{y})\in S^1\times 
{\rm I}\negthinspace {\rm R}$ and $(x,y)\in {\rm T^2.}$

\item  Let $Diff_0^{1+\epsilon }({\rm T^2})$ be the set of $C^{1+\epsilon }$
(for some $\epsilon >0$) area preserving diffeomorphisms of the torus
homotopic to the identity and let $Diff_0^{1+\epsilon }({\rm I\negthinspace %
R^2})$ be the set of lifts of elements from $Diff_0^{1+\epsilon }({\rm T^2})$
to the plane. Maps from $Diff_0^{1+\epsilon }({\rm T^2})$ are denoted $f$
and their lifts to the plane are denoted $\widetilde{f}.$

\item  Let $Diff_k^{1+\epsilon }({\rm T^2})$ be the set of $C^{1+\epsilon }$
(for some $\epsilon >0$ and some integer $k\neq 0$) area preserving
diffeomorphisms of the torus homotopic to a Dehn twist $(x,y)\longrightarrow
(x+ky$ mod$1,y$ mod$1)$ and let $Diff_k^{1+\epsilon }(S^1\times {\rm I}%
\negthinspace {\rm R})$ and $Diff_k^{1+\epsilon }({\rm I\negthinspace R^2})$
be the sets of lifts of elements from $Diff_k^{1+\epsilon }({\rm T^2})$ to
the cylinder and plane. As defined above, maps from $Diff_k^{1+\epsilon }(%
{\rm T^2})$ are denoted $f$ and their lifts to the vertical cylinder and
plane are respectively denoted $\widehat{f}$ and $\widetilde{f}.$

\item  Let $p_{1,2}:{\rm I}\negthinspace {\rm R^2}\longrightarrow {\rm I}%
\negthinspace {\rm R}$ be the standard projections; $p_1(\tilde x,\tilde
y)=\tilde x$ and $p_2(\tilde x,\tilde y)=\tilde y$. 

\item  Given $f\in Diff_0^{1+\epsilon }({\rm T^2})$ and a lift $\widetilde{f}%
\in Diff_0^{1+\epsilon }({\rm I}\negthinspace {\rm R^2}),$ the so called
rotation set of $\widetilde{f},$ $\rho (\widetilde{f}),$ can be defined as
follows (see \cite{misiu}): 
\begin{equation}
\label{rotsetident}\rho (\widetilde{f})=\bigcap_{{\ 
\begin{array}{c}
i\geq 1 \\ 
\end{array}
}}\overline{\bigcup_{{\ 
\begin{array}{c}
n\geq i \\ 
\end{array}
}}\left\{ \frac{\widetilde{f}^n(\widetilde{z})-\widetilde{z}}n:\widetilde{z}%
\in {\rm I}\negthinspace {\rm R^2}\right\} }
\end{equation}

This set is a compact convex subset of ${\rm I\negthinspace R^2}$ and it was
proved in \cite{franksrat} and \cite{misiu} that all points in its interior
are realized by compact $f$-invariant subsets of ${\rm T^2,}$ which are
periodic orbits in the rational case. We say that some vector $\rho \in \rho
(\widetilde{f})$ is realized by a compact $f$-invariant subset $K\subset 
{\rm T^2}$ if and only if, for all $z\in K$ and any $\widetilde{z}\in
p^{-1}(z)$ we have 
\begin{equation}
\label{deffrotvect}\stackunder{n\rightarrow \infty }{\lim }\frac{\widetilde{f%
}^n(\widetilde{z})-\widetilde{z}}n=\rho .
\end{equation}
The above limit, when it exists, is know as the rotation vector of the point 
$z.$ So when we say that some point $\rho $ in $\rho (\widetilde{f})$ is
realized by a compact $f$-invariant subset of ${\rm T^2}$ we mean that the
rotation vector exists at all points of this set and is equal to $\rho .$

\item  Given $f\in Diff_k^{1+\epsilon }({\rm T^2})$ and a lift $\widetilde{f}%
\in Diff_k^{1+\epsilon }({\rm I}\negthinspace {\rm R^2}),$ the so called
vertical rotation set of $\widetilde{f},$ $\rho _V(\widetilde{f}),$ can be
defined as follows, see \cite{doeff}, \cite{eu1} and \cite{eu4}: 
\begin{equation}
\label{rotsetdehn}\rho _V(\widetilde{f})=\bigcap_{{\ 
\begin{array}{c}
i\geq 1 \\ 
\end{array}
}}\overline{\bigcup_{{\ 
\begin{array}{c}
n\geq i \\ 
\end{array}
}}\left\{ \frac{p_2\circ \widetilde{f}^n(\widetilde{z})-p_2(\widetilde{z})}n:
\widetilde{z}\in {\rm I}\negthinspace {\rm R^2}\right\} }
\end{equation}


This set is a closed interval (maybe a single point, but never empty) and it
was proved in \cite{doeff}, \cite{eu1} and \cite{eu4} that all numbers in
its interior are realized by compact $f$-invariant subsets of ${\rm T^2,}$
which are periodic orbits in the rational case, in exactly the same sense as
above. The only difference is that to define the vertical rotation number of
a point $z\in {\rm T^2,}$ we consider the $p_2$-projection of the expression
(\ref{deffrotvect}). In other words we say that the vertical rotation number
of $z$ is equal to some number $\omega $ if for any $\widetilde{z}\in
p^{-1}(z),$ 
$$
\stackunder{n\rightarrow \infty }{\lim }\frac{p_2\circ \widetilde{f}^n(
\widetilde{z})-p_2(\widetilde{z})}n=\omega . 
$$

\item  A connected simply connected open subset $D$ of the torus is called
an open disk. Note that in this case, for any connected component $
\widetilde{D}$ of $p^{-1}(D)$ and any pair of integers $(a,b)\neq (0,0),$ we
have%
$$
\widetilde{D}\cap (\widetilde{D}+(a,b))=\emptyset \text{ and }p^{-1}(D)=%
\stackunder{i,j\in integers}{\cup }\widetilde{D}+(i,j).\text{ } 
$$

\item  We say that an open disk $D\subset {\rm T^2}$ is unbounded if a
connected component $\widetilde{D}$ of $p^{-1}(D)$ is unbounded. Clearly, if 
$\widetilde{D}^{\prime }$ is another connected component of $p^{-1}(D),$ as
there exists a pair of integers $(a,b)$ such that $\widetilde{D}=\widetilde{D%
}^{\prime }+(a,b),$ all connected components of $p^{-1}(D)$ are unbounded.

\item  For any hyperbolic periodic point $P$, we say that its unstable
manifold $W^u(P)$ (or stable manifold $W^s(P)$) has a topologically
transverse intersection with some closed connected set $K$ if and only if,
for some compact connected piece $\lambda $ of some branch of $W^u(P)$ (or $%
W^s(P)$), there exists a rectangle $R$ such that $R\backslash \lambda $ has
exactly two connected components, $R_{left}$ and $R_{right}$ and there is a
connected component of $K\cap R_{left}$ which intersects $\lambda $ and
another side of $R_{left}.$ An analogous condition is assumed for $%
R_{right}, $ namely there is a connected component of $K\cap R_{right}$
which intersects $\lambda $ and another side of $R_{right}.$ By rectangle we
mean that $R$ is a topological closed disk and the boundary of $R$ is a
simple closed curve divided into four $C^1$ simple arcs (the sides): two are
very small and $C^1$-transversal to $\lambda $ and the other two are $C^1$%
-close to $\lambda ,$ see figure 1. We use the following notation: $W^u(P)%
\pitchfork
K.$ It is easy to see, using Hartman-Grobman theorem, that if $P$ is a
hyperbolic fixed point for some map $f$ and $W^u(P)\pitchfork K,$ then given 
$\epsilon >0$ and a compact piece $\kappa $ of $W^s(P),$ there exists $%
n(\epsilon ,\kappa )>0$ such that the $\epsilon $-neighborhood of $f^{-n}(K)$
contains $\kappa $ for all $n>n(\epsilon ,\kappa ).$ So, if $Q$ is another
hyperbolic fixed point for $f,$ we say that $W^u(P)\pitchfork W^s(Q)$ when
there exists a closed connected piece $\alpha $ of $W^s(Q)$ such that $W^u(P)%
\pitchfork \alpha .$ In this case $\overline{W^u(P)}\supset \overline{W^u(Q)}%
,$ which implies that if $W^u(Q)\pitchfork K$ for some closed connected set $%
K,$ then $W^u(P)\pitchfork K$ (this is more or less a $C^0$ version of the $%
\lambda $-lemma and some of its consequences).
\end{enumerate}

Now we present our main results.

\begin{lemma}
\label{rothorse}: Suppose $f$ belongs to $Diff_0^{1+\epsilon }({\rm T^2})$
and $(0,0)\in int(\rho (\widetilde{f}))$ or $f$ belongs to $%
Diff_k^{1+\epsilon }({\rm T^2})$ and $0\in int(\rho _V(\widetilde{f})).$
Then there exists $\widetilde{Q}\in {\rm I}\negthinspace 
{\rm R^2,}$ which is a hyperbolic periodic point for $\widetilde{f}$ such
that for some pair of integers $(a,b)\neq (0,0),$ $a,b$ coprimes, $W^u( 
\widetilde{Q})\pitchfork W^s(\widetilde{Q}+(a,b))$ $($note that $W^s( 
\widetilde{Q}+(a,b))=W^s(\widetilde{Q})+(a,b)).$ In particular, if $Q=p( 
\widetilde{Q})$ then $W^s(Q)\cup Q\cup W^u(Q)$ contains a homotopically
non-trivial simple closed curve in ${\rm T^2.}$
\end{lemma}

\vskip 0.2truecm

%
%

This lemma is easier in case $f$ is transitive, see lemma \ref{interhet}. In
the general case, we have to work with pseudo-Anosov maps isotopic to $f$
relative to certain finite invariant sets and apply Handel's shadowing \cite
{handel} and \cite{boyland}, and other technical results on Pesin theory 
\cite{andre}.

\begin{lemma}
\label{perioddisk}: Suppose $f$ belongs to $Diff_0^{1+\epsilon }({\rm T^2})$
and $int(\rho (\widetilde{f}))$ is not empty or $f$ belongs to $%
Diff_k^{1+\epsilon }({\rm T^2})$ and $int(\rho _V(\widetilde{f}))$ is not
empty. If $f$ is transitive, then $f$ can not have a periodic open disk. In
the general case, there exists $M=M(f)>0$ such that if $D\subset {\rm T^2}$
is a $f$-periodic open disk, then for any connected component $\widetilde{D}$
of $p^{-1}(D),$ $diam(\widetilde{D})<M.$
\end{lemma}

\vskip 0.2truecm

%
%
%

So if the rotation set has interior, a diffeomorphism of the torus homotopic
to the identity or to a Dehn twist can only have bounded periodic open disks
and the bound in their diameters (from above) is uniform.

\begin{theorem}
\label{densevarident}: Suppose $f$ belongs to $Diff_0^{1+\epsilon }({\rm T^2}%
),$ it is transitive and $\left( \frac pq,\frac rq\right) \in int(\rho ( 
\widetilde{f})).$ Then, $\widetilde{f}^q(\bullet )-(p,r)$ has a hyperbolic
periodic point $\widetilde{Q}$ such that $\overline{W^u(\widetilde{Q})}= 
\overline{W^s(\widetilde{Q})}={\rm I}\negthinspace {\rm R^2.}$
\end{theorem}

\vskip 0.2truecm

\begin{theorem}
\label{densevardehn}: Suppose $f$ belongs to $Diff_k^{1+\epsilon }({\rm T^2}%
),$ it is transitive, $\frac pq\in int(\rho _V(\widetilde{f}))$ and $s$ is
any integer number. Then, $\widetilde{f}^q(\bullet )-(s,p)$ has a hyperbolic
periodic point $\widetilde{Q}$ such that $\overline{W^u(\widetilde{Q})}= 
\overline{W^s(\widetilde{Q})}={\rm I}\negthinspace {\rm R^2.}$
\end{theorem}

\vskip 0.2truecm

\begin{corollary}
\label{topmix}: Suppose $f$ belongs to $Diff_0^{1+\epsilon }({\rm T^2})$ and 
$(0,0)\in int(\rho (\widetilde{f}))$ or $f$ belongs to $Diff_k^{1+\epsilon }(%
{\rm T^2})$ and $0\in int(\rho _V(\widetilde{f})).$ If $f$ is transitive,
then $\widetilde{f}$ is topologically mixing.
\end{corollary}

\vskip 0.2truecm

%
%
%

The next result is a version of theorems \ref{densevarident} and \ref
{densevardehn} to the general case.

\begin{theorem}
\label{casogeral}: Suppose $f$ belongs to $Diff_0^{1+\epsilon }({\rm T^2})$
and $(0,0)\in int(\rho (\widetilde{f}))$ or $f$ belongs to $%
Diff_k^{1+\epsilon }({\rm T^2})$ and $0\in int(\rho _V(\widetilde{f})).$
Then, $\widetilde{f}$ has a hyperbolic periodic point $\widetilde{Q}$ such
that for any pair of integers $(a,b),$ $W^u(\widetilde{Q})\pitchfork W^s( 
\widetilde{Q}+(a,b)),$ so $\overline{W^u(\widetilde{Q})}=\overline{W^u( 
\widetilde{Q})}+(a,b)$ and $\overline{W^s(\widetilde{Q})}=\overline{W^s( 
\widetilde{Q})}+(a,b).$ Moreover, $\overline{W^u(\widetilde{Q})}=\overline{%
W^s(\widetilde{Q})},$ this set is $\widetilde{f}$-invariant ($\widetilde{f}( 
\overline{W^u(\widetilde{Q})})=\overline{W^u(\widetilde{Q})}$) and all
connected components of the complement of $\overline{W^u(\widetilde{Q})}$
are open disks with diameter uniformly bounded from above, $f$-periodic when
projected to the torus.
\end{theorem}

\vskip 0.2truecm

%
%

\begin{corollary}
\label{wandering}: Suppose $f$ belongs to $Diff_0^{1+\epsilon }({\rm T^2})$
and $(0,0)\in int(\rho (\widetilde{f}))$ or $f$ belongs to $%
Diff_k^{1+\epsilon }({\rm T^2})$ and $0\in int(\rho _V(\widetilde{f})).$ If $
\widetilde{D}\subset {\rm I}\negthinspace {\rm R^2}$ is a wandering open
disk, then for all integers $n,$ $diam(\widetilde{f}^n(\widetilde{D}))$ is
uniformly bounded from above and there exists a $f$-periodic open disk $%
D_{+}\subset {\rm T^2}$ such that $D_{+}\supset D=p(\widetilde{D}).$
\end{corollary}

\vskip 0.2truecm

Let us consider now the following subset of ${\rm T^2,}$ 
$$
R.I.(f)\stackrel{def.}{=}p(\overline{W^u(\widetilde{Q})})=\overline{W^u(p( 
\widetilde{Q}))}, 
$$
where the point $\widetilde{Q}$ appears in theorem \ref{casogeral}. This set
is called region of instability of $f.$ It is clearly compact, connected and 
$f$-invariant. And from theorem \ref{casogeral} we obtain the following:

\begin{corollary}
\label{reginst}: If $f$ belongs to $Diff_0^{1+\epsilon }({\rm T^2})$ and $%
\{(\frac pq,\frac rq),(0,0)\}\in int(\rho (\widetilde{f}))$ then $\widetilde{%
f}^q(\bullet )-(p,r)$ has a hyperbolic periodic point $\widetilde{P}$ such
that for any pair of integers $(a,b),$ $W^u(\widetilde{P})\pitchfork W^s( 
\widetilde{P}+(a,b))$ and $R.I.(f)=p(\overline{W^u(\widetilde{P})}).$ In
case $f$ belongs to $Diff_k^{1+\epsilon }({\rm T^2})$ and $\{\frac pq,0\}\in
int(\rho _V(\widetilde{f})),$ then, for all integers $s,$ $\widetilde{f}%
^q(\bullet )-(s,p)$ has a hyperbolic periodic point $\widetilde{P}$ such
that for any pair of integers $(a,b),$ $W^u(\widetilde{P})\pitchfork W^s( 
\widetilde{P}+(a,b))$ and $R.I.(f)=p(\overline{W^u(\widetilde{P})}).$
Moreover, if $\widetilde{Q}$ is the $\widetilde{f}$-periodic point that
appears in the statement of theorem \ref{casogeral}, then $W^u(p(\widetilde{P%
}))\pitchfork W^s(p(\widetilde{Q}))$ and $W^u(p(\widetilde{Q}))\pitchfork %
W^s(p(\widetilde{P})).$
\end{corollary}

We believe this is an interesting description of the region of instability
of an area-preserving diffeomorphism of the torus. Moreover, a result
analogous to corollary \ref{topmix} holds. The map $\widetilde{f}$
restricted to $\overline{W^u(\widetilde{Q})}=\overline{W^s(\widetilde{Q})}$
is topologically mixing, that is for every $\widetilde{z}_1,\widetilde{z}%
_2\in \overline{W^u(\widetilde{Q})}$ and $\epsilon >0,$ there exists $N>0$
such that 
$$
\widetilde{f}^n(B_\epsilon (\widetilde{z}_1))\cap B_\epsilon (\widetilde{z}%
_2)\neq \emptyset ,\forall n\geq N. 
$$
Clearly, many interesting and fundamental questions remain open: Is $f$
transitive when restricted to $R.I.(f)?$ Does $f$ has a dense set of
periodic points in $R.I.(f)?$ What about $\widetilde{f}$ restricted to $%
p^{-1}(R.I.(f))?$ Is it transitive? Does it has a dense set of periodic
points?

%
%

One consequence of the above results, that may be useful in many situations
will be presented. For this, we need more definitions.

In the homotopic to the identity case, given a vector $(\cos (\theta ),\sin
(\theta )),$ we define 
$$
B_\theta =\{\widetilde{z}\in {\rm I}\negthinspace {\rm R^2:}\text{ }%
\left\langle \widetilde{f}^n(\widetilde{z}),(\cos (\theta ),\sin (\theta
))\right\rangle \geq 0\text{ for all integers }n\geq 0\} 
$$
and let $B_\theta ^\infty $ be the union of all unbounded connected
components of $B_\theta .$ In case $f$ is homotopic to a Dehn twist, it only
makes sense to define the sets

$$
B_{S(\text{or }N)}=\{\widetilde{z}\in {\rm I}\negthinspace {\rm R^2:}\text{ }%
p_2\circ \widetilde{f}^n(\widetilde{z})\leq (\geq )0\text{ for all integers }%
n\geq 0\} 
$$
and let $B_{S(\text{or }N)}^\infty $ be the union of all unbounded connected
components of $B_{S(\text{or }N)}.$

In lemma 2 of \cite{eufabraulio} we proved that, if $(0,0)\in \rho ( 
\widetilde{f})$ or $0\in \rho _V(\widetilde{f}),$ then $B_\theta ^\infty ,$ $%
B_S^\infty $ and $B_N^\infty $ are non-empty, closed subsets of the plane,
positively invariant under $\widetilde{f}.$ It is easy to see that their
omega-limits satisfy the following (because $\widetilde{f}(B_\theta ^\infty
)\subset B_\theta ^\infty $ and $\widetilde{f}(B_{S(N)}^\infty )\subset
B_{S(N)}^\infty ):$ 
$$
\begin{array}{c}
\\ 
\omega (B_\theta ^\infty )= 
\stackunder{i\geq 0}{\cap }\widetilde{f}^i(B_\theta ^\infty )=\stackunder{%
i\in integers}{\cap }\widetilde{f}^i(B_\theta ^\infty ) \\  \\ 
\omega (B_{S(N)}^\infty )=\stackunder{i\geq 0}{\cap }\widetilde{f}%
^i(B_{S(N)}^\infty )=\stackunder{i\in integers}{\cap }\widetilde{f}%
^i(B_{S(N)}^\infty ) 
\end{array}
$$

And we know that if $\omega (B_\theta ^\infty )=\emptyset ,$ then all points 
$\widetilde{z}\in B_\theta ^\infty $ satisfy (see lemma 10 of \cite{eufa1})%
$$
\left\langle \left( \frac{\widetilde{f}^n(\widetilde{z})-\widetilde{z}}%
n\right) ,(\cos (\theta ),\sin (\theta ))\right\rangle >c_\theta >0, 
$$
for some constant $c_\theta $ and all $n>0$ suff. large. Analogously, if $%
\omega (B_{S(N)}^\infty )=\emptyset ,$ then all points $\widetilde{z}\in
B_{S(N)}^\infty $ satisfy%
$$
\frac{p_2\circ \widetilde{f}^n(\widetilde{z})-p_2(\widetilde{z})}n<c_S<0 
\text{ (resp. }>c_N>0), 
$$
for some constant $c_S$ ($c_N$) and all $n>0$ suff. large.

\begin{theorem}
\label{omegaident}: Suppose $f$ belongs to $Diff_0^{1+\epsilon }({\rm T^2})$
and $(0,0)\in int(\rho (\widetilde{f}))$ or $f$ belongs to $%
Diff_k^{1+\epsilon }({\rm T^2})$ and $0\in int(\rho _V(\widetilde{f})).$
Then, for all $\theta \in [0,2\pi ],$ $\omega (B_\theta ^\infty )=\emptyset $
and $\omega (B_S^\infty )=\omega (B_N^\infty )=\emptyset .$
\end{theorem}

\vskip 0.2truecm

%
%

A simple corollary of the above result, lemma \ref{perioddisk} and the ideas
in the proof of theorem \ref{casogeral} is the following:

\begin{corollary}
\label{sobreb}: Suppose $f$ belongs to $Diff_0^{1+\epsilon }({\rm T^2})$ and 
$(0,0)\in int(\rho (\widetilde{f}))$ or $f$ belongs to $Diff_k^{1+\epsilon }(%
{\rm T^2})$ and $0\in int(\rho _V(\widetilde{f})).$ If $f$ is transitive,
then $\overline{p(B_S^\infty )}=\overline{p(B_N^\infty )}=\overline{%
p(B_\theta ^\infty )}={\rm T^2,}$ for all $\theta \in [0,2\pi ].$ And for a
general $f,$ any of the following sets $\left( \overline{p(B_S^\infty )}%
\right) ^c,$ $\left( \overline{p(B_N^\infty )}\right) ^c$ and $\left( 
\overline{p(B_\theta ^\infty )}\right) ^c$ is the union of $f$-periodic open
disks, with diameters uniformly bounded from above when lifted to the plane.
\end{corollary}

The hypothesis of zero being an interior point of the rotation set is
essential in theorem \ref{omegaident}, as the following example shows.
Consider the standard map $S_{M,k_0}$ for a parameter $k_0>0$ which has a
rotational invariant curve $\gamma .$ Now, if $\epsilon >0$ is sufficiently
small, the map $S_{M,k_0}^\epsilon :{\rm T^2\rightarrow T^2,}$ 
$$
S_{M,k_0}^\epsilon (x,y)=(x+y+k_0\sin (2\pi x)\text{ }mod1,y+k_0\sin (2\pi
x)+\epsilon \text{ }mod1) 
$$
still has fixed points of zero vertical rotation number (because $S_{M,k_0}$
has a hyperbolic fixed point with zero vertical rotation number)\ and, as
the vertical rotation number of Lebesgue measure 
$$
\int_{{\rm T^2}}(k_0\sin (2\pi x)+\epsilon )dLeb=\epsilon >0, 
$$
we get that $\rho _V(\widetilde{S}_{M,k_0}^\epsilon )$ is an interval which
has zero as the left extreme. Note that $\widetilde{S}_{M,k_0}^\epsilon :%
{\rm I}\negthinspace {\rm R^2\rightarrow I}\negthinspace {\rm R^2}$ is given
by%
$$
\widetilde{S}_{M,k_0}^\epsilon (\widetilde{x},\widetilde{y})=(\widetilde{x}+ 
\widetilde{y}+k_0\sin (2\pi \widetilde{x}),\widetilde{y}+k_0\sin (2\pi 
\widetilde{x})+\epsilon ). 
$$
This happens because the extremes of the vertical rotation interval of $
\widetilde{S}_{M,k_0}^\epsilon $ are continuous non-decreasing functions of $%
\epsilon >0$ by theorem 10 and lemma 2 of \cite{eucont}. 

Now we will show that $\omega (B_S^\infty )$ and $\omega (B_N^\infty )$ of $
\widetilde{S}_{M,k_0}^\epsilon$ are both non-empty. Birkhoff's invariant
curve theorem implies that $\gamma $ projects injectively on the horizontal
direction, so if we consider the cylinder diffeomorphism induced by $
\widetilde{S}_{M,k_0}^\epsilon ,$ denoted $\widehat{S}_{M,k_0}^\epsilon
:S^1\times {\rm I}\negthinspace {\rm R\rightarrow }S^1\times {\rm I}%
\negthinspace {\rm R,\ }$we get that the closure of the region below $
\widehat{\gamma }$ (denoted $\widehat{\gamma }^{-}$) is invariant under $( 
\widehat{S}_{M,k_0}^\epsilon )^{-1}$ and the closure of the region above $
\widehat{\gamma }$ (denoted $\widehat{\gamma }^{+}$) is invariant under $
\widehat{S}_{M,k_0}^\epsilon .$ Note that $\widehat{\gamma }$ is just a lift
of $\gamma $ to the cylinder. Moreover, as $\widehat{S}_{M,k_0}^\epsilon $
has fixed points above and below $\widehat{\gamma },$ the sets 
$$
\stackunder{i\geq 0}{\cap }(\widehat{S}_{M,k_0}^\epsilon )^i(\widehat{\gamma 
}^{+})\text{ and }\stackunder{i\geq 0}{\cap }(\widehat{S}_{M,k_0}^\epsilon
)^{-i}(\widehat{\gamma }^{-})\text{ } 
$$
are closed, non-empty, $\widehat{S}_{M,k_0}^\epsilon $-invariant and their
connected components are unbounded. Clearly if we consider their inverse
images under the projection $\pi :{\rm I}\negthinspace {\rm R^2\rightarrow }%
S^1\times {\rm I}\negthinspace {\rm R,}$ and take adequate vertical
translations, we get that both $\omega (B_S^\infty )$ and $\omega
(B_N^\infty )$ are non-empty.

This paper is organized as follows. In the second section we present some
background results we use, with references and in the third section we prove
our results. Just to avoid confusions, let us state that everytime we say a
point $Q\in {\rm T^2}$ is $f$-periodic of period $m$ and hyperbolic, we not
only mean that $f^m(Q)=Q,$ but also that $Df^m(P)$ has positive eigenvalues,
which implies that all branches of the stable and unstable manifolds at $Q$
are $f^m$-invariant.

\section{Ideas involved in the proofs}

The main tools used in our proofs come from two different theories. We will
try to give a superficial description on each of them below.

\subsection{Nielsen-Thurston theory of classification of homeomorphisms of
surfaces and some results on subshifts of finite type}

The following is a brief summary of this powerful theory. For more
information and proofs see \cite{T}, \cite{FLP} and \cite{HT}.

Let $M$ be a compact, connected oriented surface possibility with boundary,
and $f:M\rightarrow M$ be a homeomorphism. Two homeomorphisms are said to be
isotopic if they are homotopic via homeomorphisms. In fact, for closed
orientable surfaces, all homotopic pairs of homeomorphisms are isotopic \cite
{E1}.

There are two basic types of homeomorphisms which appear in the
Nielsen-Thurston classification\ :\ the finite order homeomorphisms and the
pseudo-Anosov ones.

A homeomorphism $f$ is said to be of finite order if $f^n=id$ for some $n\in 
{\rm I\negthinspace N.}$ The least such $n$ is called the order of $f.$
Finite order homeomorphisms have zero topological entropy.

A homeomorphism $f$ is said to be pseudo-Anosov if there is a real number $%
\lambda >1$ and a pair of transverse measured foliations ${\cal F}^u$ and $%
{\cal F}^s$ such that $f({\cal F}^s)=\lambda ^{-1}{\cal F}^s$ and $f({\cal F}%
^u)=\lambda {\cal F}^u.$ Pseudo-Anosov homeomorphisms preserve area, are
topologically transitive, have positive topological entropy, and have Markov
partitions \cite{FLP}.

A homeomorphism $f$ is said to be reducible by a system 
$$
C=\stackrel{n}{\stackunder{i=1}{\cup }}C_i 
$$

of disjoint simple closed curves $C_1,...,C_n$ (called reducing curves) if

(1) $\forall i,$ $C_i$ is not homotopic to a point, nor to a component of $%
\partial M,$

(2) $\forall i\neq j,$ $C_i$ is not homotopic to $C_j,$

(3) $C$ is invariant under $f.$

\begin{theorem}
\label{Thu}\ :\ If the Euler characteristic $\chi (M)<0,$ then every
homeomorphism $f:M\rightarrow M$ is isotopic to a homeomorphism $\phi
:M\rightarrow M$ such that either

(a) $\phi $ is of finite order,

(b) $\phi $ is pseudo-Anosov, or

(c) $\phi $ is reducible by a system of curves $C.$
\end{theorem}

Maps $\phi $ as in theorem \ref{Thu} are called Thurston canonical forms for 
$f.$

Now two applications of theorem \ref{Thu} that are important to us will be
presented. The first is due to Llibre and Mackay \cite{llibre} and Franks 
\cite{franksrat}, and the second is due to Doeff \cite{doeff} and myself 
\cite{eu1}:

\begin{theorem}
\label{llibremack}: If $f$ is a homeomorphism of the torus homotopic to the
identity and $\widetilde{f}$ is a lift of $f$ to the plane such that $\rho ( 
\widetilde{f})$ has interior, then for any three non-collinear rational
vectors $\rho _1,\rho _2,\rho _3\in int(\rho (\widetilde{f})),$ there are
periodic orbits $Q_1,Q_2$ and $Q_3$ that realize these rotation vectors such
that $f\mid _{{\rm T^2}\backslash \{Q_1\cup Q_2\cup Q_3\}}$ is isotopic to a
pseudo-Anosov homeomorphism $\phi $ of ${\rm T^2}\backslash \{Q_1\cup
Q_2\cup Q_3\}.$
\end{theorem}

In this case we say that $f$ is isotopic to $\phi $ relative to $Q_1\cup
Q_2\cup Q_3.$ It means that the isotopy acts on the set $Q_1\cup Q_2\cup Q_3$
exactly in the same way as $f$ does. The map $\phi $$:{\rm T^2\rightarrow T^2%
}$ is then said to be pseudo-Anosov relative to the finite invariant set $%
Q_1\cup Q_2\cup Q_3$ because it satisfies all of the properties of a
pseudo-Anosov homeomorphism except that the associated stable and unstable
foliations may have 1-prong singularities at points in $Q_1\cup Q_2\cup Q_3.$

\begin{theorem}
\label{meuedoeff}: If $f$ is a homeomorphism of the torus homotopic to a
Dehn twist and $\widetilde{f}$ is a lift of $f$ to the plane such that $\rho
_V(\widetilde{f})$ has interior, then for any two different rationals $\rho
_1,\rho _2\in int(\rho _V(\widetilde{f})),$ there are periodic orbits $Q_1$
and $Q_2$ that realize these rotation numbers such that $f$ is isotopic
relative to $Q_1\cup Q_2$ to a homeomorphism $\phi $ of ${\rm T^2,}$ which
is pseudo-Anosov relative to $Q_1\cup Q_2.$
\end{theorem}

In theorem \ref{llibremack}, given a rational vector $\rho =(\frac pq,\frac
rq)$ such that the integers $p,r,q$ have no common factors, the chosen
periodic orbit $Q$ which realizes this rotation vector must have period $q.$
The existence of such an orbit follows from the main theorem in \cite
{franksrat}. An analogous remark holds for theorem \ref{meuedoeff}, namely
given a rational $\frac pq$ $\in int(\rho _V(\widetilde{f}))$ such that the
integers $p,q$ have no common factors, the chosen periodic orbit $Q$ which
realizes this vertical rotation number must also have period $q.$ The
existence of such an orbit follows from the main theorems in \cite{doeff}
and \cite{eu4}.

Using the notation from theorem \ref{Thu}, our next result tells us that if $%
\phi $ is pseudo-Anosov, then the complicated dynamics of $\phi $ is in some
sense shared by $f.$ To be more precise, we will present the following
analog of Handel's global shadowing \cite{handel}, see theorem 3.3.1 of \cite
{kwapisz}, which is more or less taken from theorem 3.2 of \cite{boyland},
stated in a version adequate to our needs:

\begin{theorem}
\label{handel}: If $f:{\rm T^2\rightarrow T^2}$ is a diffeomorphism
homotopic to the identity or to a Dehn twist, $A$ is a finite $f$-invariant
set and $f$ is isotopic rel $A$ to some map $\phi :{\rm T^2\rightarrow T^2}$
which is pseudo-Anosov rel $A,$ then there exists a compact $f$-invariant
set $W\subset {\rm T^2}$ and a continuous surjection $s:W\rightarrow {\rm T^2%
}$ that is homotopic to the inclusion map $i:W\rightarrow {\rm T^2}$ and
semi-conjugates $f\mid _W$ to $\phi ,$ that is, $s\circ f\mid _W=\phi \circ
s.$
\end{theorem}

{\bf Remark: } The set $W,$ which is the domain of $s,$ is the closure of
another set denoted $W^{\prime }$ that satisfies the following. The map $s$
restricted to $W^{\prime }$ is one to one and $s(W^{\prime })=\{periodic$ $%
points$ $of$ $\phi \},$ which is a dense subset of ${\rm T^2,}$ see the
proof of theorem 3.2 of \cite{boyland}. In this way, given a compact subset $%
M$ of $W,$ if $s(M)={\rm T^2,}$ then $M=W.$

As $s:W\rightarrow {\rm T^2}$ is homotopic to the inclusion map $%
i:W\rightarrow {\rm T^2,}$ we get that $s$ has a lift $\widetilde{s}%
:p^{-1}(W)\rightarrow {\rm I}\negthinspace {\rm R^2}$ such that%
$$
\widetilde{s}\circ \widetilde{f}\mid _{p^{-1}(W)}=\widetilde{\phi }\circ 
\widetilde{s}\text{ and }\stackunder{\widetilde{z}\in p^{-1}(W)}{\sup }%
\left\| \widetilde{s}(\widetilde{z})-\widetilde{z}\right\| <\infty , 
$$
where $\widetilde{f}$ and $\widetilde{\phi }$ are lifts that are
equivariantly homotopic rel $p^{-1}(A).$

The next results we present are related to theorems \ref{llibremack} and \ref
{meuedoeff} in the following sense: Assume $\phi :{\rm T^2\rightarrow T^2}$
is a homeomorphism of the torus, either homotopic to the identity or to a
Dehn twist, which is pseudo-Anosov relative to some finite $\phi $-invariant
set $Q.$ Also suppose that its rotation set, or vertical rotation set
(depending on whether $\phi $ is homotopic to the identity or to a Dehn
twists) has $0$ as an interior point. The first result appears in \cite{FLP}
and \cite{kwapisz}:

\begin{theorem}
\label{flp}: There exists a Markov partition ${\cal R}=\{R_1,...,R_N\}$ for $%
\phi .$ If $G$ is a graph with the set of vertices $\{1,...,N\}$ that has an
edge from $i$ to $j$ whenever $\phi (R_i)\cap int(R_j)\neq \emptyset ,$ then
the subshift of finite type $(\Lambda ,\sigma )$ associated to $G$ is mixing
and factors onto $\phi .$ More precisely, for any $x=(x_i)_{i\in
integers}\in \Lambda $ the intersection $\stackunder{i\in integers}{\cap }%
\phi ^{-i}(R_{x_i})$ consists of a single point, denoted $h(x)$ and the map $%
h:\Lambda \rightarrow {\rm T^2}$ has the following properties:

\begin{enumerate}
\item  $h\circ \sigma =\phi \circ h;$

\item  $h$ is continuous and surjective;

\item  $h$ is finite to one;

\item  $h$ is one to one on a topologically residual subset consisting of
all points whose full orbits never hit the boundary of the Markov partition;
\end{enumerate}
\end{theorem}

The subshift of finite type that appears in the above theorem can be defined
as follows: Let $A$ be a $N\times N$ matrix whose entries satisfy $a_{ij}=1$
if and only if $\phi (R_i)\cap int(R_j)\neq \emptyset $ and $a_{ij}=0$
otherwise. Following the theorem notation, $\Lambda $ is the space of
bi-infinite sequences $x=(...,x_{-1},x_0,x_1,...)\in \{1,...,N\}^{{\rm Z%
\negthinspace
\negthinspace Z}}$ such that $a_{x_lx_{l+1}}=1$ for all integers $l$ and $%
\sigma :\Lambda \rightarrow \Lambda $ is just the shift transformation $%
\sigma (x)_n=(x)_{n+1}.$

Now, assume $\psi :\Lambda \rightarrow {\rm I}\negthinspace {\rm R^2}$ is a
function depending only on $x_0$ and $x_1$ for which there exists an
Euclidean triangle $\Delta \subset {\rm I}\negthinspace {\rm R^2}$ such that 
$$
\Delta \subset \{\int_\Lambda \psi d\nu :\text{ }\nu \text{ is a }\sigma 
\text{-invariant Borel probability measure}\} 
$$
and for all $\rho \in int(\Delta )$ there exists a $\sigma $-invariant
ergodic Borel probability measure $\nu $ that satisfies $\int_\Lambda \psi $ 
$d\nu =\rho .$ Now we state a theorem which is contained in theorems 2.2.1
and 2.2.2 of \cite{kwapisz}:

\begin{theorem}
\label{subshift}: Under the above assumptions, for every $\rho \in
int(\Delta ),$ there exists a compact $\sigma $-invariant set $K_\rho
\subset \Lambda $ such that $h_{top}(\sigma \mid _{K_\rho })>0$ (the
topological entropy of $\sigma $ restricted to $K_\rho $ is non zero) and
for some constant $Const>0,$ all $x\in K_\rho $ and all integers $n>0$ we
have:%
$$
\left\| \stackrel{n-1}{\stackunder{j=0}{\sum }}\psi \circ \sigma
^j(x)-n.\rho \right\| \leq Const 
$$
\end{theorem}

\subsection{Pesin-Katok theory}

This is a series of amazing results, which show the existence of a certain
type of hyperbolicity everytime a diffeomorphism is sufficiently smooth (at
least $C^{1+\epsilon },$ for some $\epsilon >0$) and satisfies certain
conditions on invariant measures. Unfortunately results on this subject are
quite technical and the definitions needed are really complicated. For this
reason, we decided to present a very naive description of the theory in the
particular case of surfaces and then state a theorem.

In case of surface diffeomorphisms, positive topological entropy is the
crucial hypothesis, because of the following. By the entropy variational
principle, the so called Ruelle-Pesin inequality and the fact that ergodic
measures are extreme points of the set of Borel probability invariant
measures, when topological entropy is positive, there always exist ergodic
invariant measures $\mu $ with non-zero Lyapunov exponents, one positive and
one negative and positive metric entropy, $h_\mu (f)>0$ (see for instance 
\cite{walters} and \cite{livrokatok}). These measures are called hyperbolic
measures. Below we state a theorem adapted from Katok's work on the subject.
For proofs, see \cite{artkatok} and the supplement of \cite{livrokatok} by
Katok-Mendoza.

\begin{theorem}
\label{katokok}: Let $f$ be a $C^{1+\epsilon }$ (for some $\epsilon >0$)
diffeomorphism of a surface $M$ and suppose $\mu $ is an ergodic hyperbolic
Borel probability $f$-invariant measure with $h_\mu (f)>0$ and compact
support. Then, for any $\alpha >0$ and any $x\in supp(\mu ),$ there exists a
hyperbolic periodic point $Q\in B_\alpha (x)$ which has a transversal
homoclinic intersection and the whole orbit of $Q$ is contained in the $%
\alpha $-neighborhood of $supp(\mu ).$
\end{theorem}

%
%
%

\section{Proofs}

In this section we will prove our main results. With this purpose we will
prove several auxiliary propositions and lemmas.

\subsection{Some preliminary results}

The first one is a by product of some results in Kwapisz thesis \cite
{kwapisz} and Katok's work on Pesin theory, see theorem \ref{katokok}.

\begin{lemma}
\label{hiperident}: Suppose $f$ belongs to $Diff_0^{1+\epsilon }({\rm T^2}).$
Then for any rational vector $\left( \frac pq,\frac rq\right) \in int(\rho ( 
\widetilde{f}))$ (with $q>0$)$,$ $\widetilde{f}^q(\bullet )-(p,r)$ has a
hyperbolic periodic point $\widetilde{Q}\in {\rm I}\negthinspace {\rm R^2}$
and $W^u(\widetilde{Q})$ has a transverse intersection with $W^s(\widetilde{Q%
}).$
\end{lemma}

{\it Proof:}

Given $\left( \frac pq,\frac rq\right) \in int(\rho (\widetilde{f})),$ let $
\widetilde{g}\stackrel{def.}{=}\widetilde{f}^q(\bullet )-(p,r).$ It is easy
to see that $(0,0)\in int(\rho (\widetilde{g}))$ and a periodic point for $
\widetilde{g}$ corresponds to a periodic point for $f$ with rotation vector
(for $\widetilde{f}$) equal $\left( \frac pq,\frac rq\right) .$ So without
loss of generality, we can suppose that $\left( \frac pq,\frac rq\right)
=(0,0)\in int(\rho (\widetilde{f})).$

Choose three periodic orbits $Q_1,Q_2$ and $Q_3$ as in theorem \ref
{llibremack}, such that their rotation vectors form a triangle $\Delta $
that contains $(0,0)$ in its interior. Theorem \ref{llibremack} tells us
that the isotopy class of $f$ relative to $Q_1\cup Q_2\cup Q_3$ contains a
pseudo-Anosov map $\phi :{\rm T^2\rightarrow T^2}$ rel $Q_1\cup Q_2\cup Q_3$%
. As we already said, pseudo-Anosov maps have very rich dynamics, for
instance theorem \ref{flp} holds.

Moreover, as $\phi :{\rm T^2\rightarrow T^2}$ is a homeomorphism of the
torus homotopic to the identity, there exists a lift of $\phi $ to the
plane, denoted $\widetilde{\phi }:{\rm I}\negthinspace {\rm R^2\rightarrow I}%
\negthinspace {\rm R^2}$ which is equivariantly homotopic to $\widetilde{f}$
rel $p^{-1}(Q_1\cup Q_2\cup Q_3)$ and its rotation set $\rho (\widetilde{%
\phi })\supset \Delta .$ The rotation set of $\widetilde{\phi }$ can also be
obtained in the following way: Let $D_\phi :{\rm T^2\rightarrow I}%
\negthinspace {\rm R^2}$ be the displacement function, given by $D_\phi
(x,y)=\widetilde{\phi }(\widetilde{x},\widetilde{y})-(\widetilde{x}, 
\widetilde{y}),$ for any $(\widetilde{x},\widetilde{y})\in p^{-1}(x,y).$
Then 
\begin{equation}
\label{outradefrotset}\rho (\widetilde{\phi })=\{\int_{{\rm T^2}}D_\phi d\mu
:\text{ }\mu \text{ is a }\phi \text{-invariant Borel probability measure}%
\}. 
\end{equation}

As we already said when we presented the definition of rotation set, in \cite
{franksrat} and \cite{misiu} it is proved that every interior point $\omega $
of the rotation set is realized by a compact $f$-invariant subset of the
torus and so by Krylov-Boguliubov's theorem (see \cite{walters}) there
exists an ergodic $f$-invariant measure $\mu $ which realizes this rotation
vector, that is $\omega =\int_{{\rm T^2}}D_\phi $ $d\mu $. In \cite{misiu}
it is also shown that if $\omega $ is an extreme point of the rotation set,
which is a convex subset of the plane, then there also exists an ergodic $f$%
-invariant measure $\mu $ which realizes $\omega .$ This justifies why the
equality in (\ref{outradefrotset}) is true.

Using the notation introduced in theorem \ref{flp}, it is not very difficult
to see (see claim 3.2.1 of \cite{kwapisz}) that we can define a function $%
\psi :\Lambda \rightarrow {\rm I}\negthinspace {\rm R^2}$ (depending only on 
$x_0$ and $x_1)$ such that for any $x=(x_i)_{i\in integers}\in \Lambda $ and
any natural $n,$%
\begin{equation}
\label{blab}\left\| \stackrel{n-1}{\stackunder{j=0}{\sum }}\psi \circ \sigma
^j(x)-\stackrel{n-1}{\stackunder{j=0}{\sum }}D_\phi \circ \phi
^j(h(x))\right\| \leq 2\stackunder{1\leq i\leq N}{\max }\{diam(R_i)\}. 
\end{equation}

As we explained above, for all $\rho \in int(\Delta )$ there exists a $\phi $%
-invariant ergodic Borel probability measure $\mu $ such that $\int_{{\rm T^2%
}}D_\phi $ $d\mu =\rho .$ So if we consider the measure $\nu (\bullet )%
\stackrel{def.}{=}\mu (h(\bullet )),$ it is a $\sigma $-invariant ergodic
Borel probability measure which satisfies $\int_\Lambda \psi $ $d\nu =\rho .$
Thus if we fix $\rho =(0,0),$ we get from theorems \ref{flp}, \ref{subshift}
and expression (\ref{blab}) that there exists a compact $\widetilde{\phi }$%
-invariant set $K_\phi \subset {\rm I}\negthinspace {\rm R^2}$ such that $%
h_{top}(\widetilde{\phi }\mid _{K_\phi })>0.$ We want to have a similar
statement for the map $\widetilde{f}.$ With this purpose, we use theorem \ref
{handel} and the remark after it.

In this way, we obtain a compact $\widetilde{f}$-invariant set $K_f= 
\widetilde{s}^{-1}(K_\phi )\subset {\rm I}\negthinspace {\rm R^2}$ such that 
$h_{top}(\widetilde{f}\mid _{K_f})>0.$ As $h_{top}(\widetilde{f}\mid
_{K_f})>0,$ there exists a hyperbolic ergodic Borel probability $\widetilde{f%
}$-invariant measure $\mu ,$ with positive metric entropy, whose support is
contained in $K_f.$ So, given $\epsilon >0,$ theorem \ref{katokok} implies
that $\widetilde{f}$ has a hyperbolic periodic point $\widetilde{Q}$ with a
transversal homoclinic intersection and the whole orbit of $\widetilde{Q}$
is contained in $V_\epsilon (supp(\mu ))\subset V_\epsilon (K_f).$ $\Box $

\vskip 0.2truecm The next result is an analog of lemma \ref{hiperident} for
maps homotopic to Dehn twists.

\begin{lemma}
\label{hiperdehn}: Suppose $f$ belongs to $Diff_k^{1+\epsilon }({\rm T^2}).$
Then for any rational $\frac pq\in int(\rho _V(\widetilde{f}))$ (with $q>0$)
and any integer $s,$ $\widetilde{f}^q(\bullet )-(s,p)$ has a hyperbolic
periodic point with a transversal homoclinic intersection.
\end{lemma}

{\it Proof:}

As in the previous lemma, without loss of generality, we can suppose that $%
(q,p,s)=(1,0,0).$ Here we will present a different argument, that can also
be used to prove lemma \ref{hiperident}. The reason for this is the
following. Using the ideas from the proof of lemma \ref{hiperident} in this
context, we could not have any control on the integer $s$ that appears in
the statement of the present lemma. We just could prove that if $f\in
Diff_k^{1+\epsilon }({\rm T^2})$ and $\frac pq\in int(\rho _V(\widetilde{f}%
)),$ then $f$ has a hyperbolic periodic point with a transversal homoclinic
intersection, whose vertical rotation number is $p/q.$ The control on $s$ is
necessary, for instance, to prove corollary \ref{topmix}. So let us start
the argument.

Choose two periodic orbits $Q_1$ and $Q_2$ as in theorem \ref{meuedoeff},
such that their vertical rotation numbers, $\rho _1$ and $\rho _2$ are, one
positive and one negative. Theorem \ref{meuedoeff} tells us that $f$ is
isotopic relative to $Q_1\cup Q_2$ to a pseudo-Anosov map $\phi :{\rm %
T^2\rightarrow T^2}$ rel $Q_1\cup Q_2$. Clearly $\phi :{\rm T^2\rightarrow
T^2}$ is an area preserving homeomorphism of the torus homotopic to a Dehn
twist. So, there exists a lift of $\phi $ to the plane, denoted $\widetilde{%
\phi }:{\rm I}\negthinspace {\rm R^2\rightarrow I}\negthinspace {\rm R^2}$
which is equivariantly homotopic to $\widetilde{f}$ rel $p^{-1}(Q_1\cup Q_2)$
and its vertical rotation set satisfies $0\in ]\rho _2,\rho _1[\subset
int(\rho _V(\widetilde{\phi })).$ This implies by theorem 6 of \cite{eu4}
that $\widetilde{\phi }$ has fixed points. %
As $\phi $ is a pseudo-Anosov homeomorphism of the torus relative to some
finite set, the foliations ${\cal F}^u,{\cal F}^s$ may have a finite number
of singularities. Some of these singularities are $p$-prong singularities,
for some $p\geq 3$ and in $Q_1\cup Q_2$ the foliations may have $1$-prong
singularities, see \cite{FLP} and \cite{handelrot}. So if $P\in {\rm T^2}$
is a $\phi $-periodic point whose vertical rotation number belongs to $]\rho
_2,\rho _1[,$ the dynamics of some adequate iterate of $\phi $ near $P$ is
generated by finitely many invariant hyperbolic sectors glued together. In
each sector the dynamics is locally like the dynamics in the first quadrant
of the map $(\widetilde{x},\widetilde{y})\rightarrow (\alpha \widetilde{x}%
,\beta \widetilde{y}),$ for some real numbers $0<\beta <1<\alpha .$ The main
difference from the dynamics in a neighborhood of a hyperbolic periodic
point of a two-dimensional diffeomorphism is the fact that there may be more
than four hyperbolic sectors (when $P$ coincides with a $p$-prong
singularity, for some $p\geq 3)$, but never less because the vertical
rotation number of points in $Q_1\cup Q_2$ belongs to $\{\rho _2,\rho _1\}.$

\begin{proposition}
\label{ptohiperpA}: For any rational $\frac rn\in ]\rho _2,\rho _1[$ (with $%
n>0$) and any integer $s,$ $\widetilde{\phi }^n(\bullet )-(s,r)$ has a
hyperbolic periodic point with a transversal homoclinic intersection when
projected to the torus.
\end{proposition}

{\bf Remarks:}

\begin{enumerate}
\item  {\bf \ }When we say hyperbolic in the statement of this proposition,
we mean that the local dynamics is obtained by gluing exactly four sectors,
or equivalently the point is a regular point of the foliations.

\item  An analogous result holds in the homotopic to the identity case,
namely suppose $g$ belongs to $Diff_0^{1+\epsilon }({\rm T^2})$ and $%
(0,0)\in int(\rho (\widetilde{g})).$ Then, as we did in lemma \ref
{hiperident} we can choose three periodic orbits $P_1,P_2$ and $P_3$ such
that their rotation vectors form a triangle $\Delta $ that contains $(0,0)$
in its interior and the isotopy class of $g$ relative to $P_1\cup P_2\cup
P_3 $ contains a pseudo-Anosov map $\varphi :{\rm T^2\rightarrow T^2}$ rel $%
P_1\cup P_2\cup P_3$. The analogous version of proposition \ref{ptohiperpA}
to this situation is the following. {\it For any rational $(\frac pn,\frac
rn)\in int(\Delta ),$ $\widetilde{\varphi }^n(\bullet )-(p,r)$ has a
hyperbolic periodic point with a transversal homoclinic intersection when
projected to the torus. }The same proof works in both cases.
\end{enumerate}

{\it Proof of proposition \ref{ptohiperpA}: }

Without loss of generality, we can suppose that $(n,r,s)=(1,0,0).$ As we
already said, from theorem 6 of \cite{eu4} $\widetilde{\phi }$ has fixed
points. As $\phi $ is pseudo-Anosov relative to some finite set, the fixed
points of $\widetilde{\phi }$ project to the torus into a finite set $K_1.$
By the generalized Lefschetz index formula and the fact that $\phi $ is
homotopic to a Dehn twist, the sum of the topological indexes of $\phi $ on
these fixed points is zero. This follows from the fact that 
$$
\det \left( Id-\left( 
\begin{array}{cc}
1 & k \\ 
0 & 1 
\end{array}
\right) \right) =0, 
$$
so all Nielsen classes are non-essential for $\phi ,$ see for instance \cite
{livrokatok}, more precisely the comments after theorem 8.7.1.

But for some appropriate sufficiently large integer $m_2>0,$ the local
dynamics at points in $K_1$ implies that 
$$
\stackunder{z\in K_1}{\sum }ind(\phi ^{m_2},z)<0. 
$$
This happens because all points in $K_1$ with negative index are saddles
(maybe with more than four sectors) and points with positive index are
rotating saddles. By this we mean that the hyperbolic sectors around the
point rotate under iterations of $\phi $ until they fall on themselves. The
orientation reversing saddle $(\widetilde{x},\widetilde{y})\rightarrow (-2 
\widetilde{x},-0.3\widetilde{y})$ is an example such that for it, $m_2=2.$
So, for an adequate sufficiently large iterate of $\phi ,$ all points in $%
K_1 $ are saddles with each separatrix fixed, and thus they all have
negative index.

Let us look at all the fixed points of $\widetilde{\phi }^{m_2}.$ Clearly,
when projected to the torus, this set, denoted $K_2\supset K_1,$ is also
finite and not equal to $K_1$ because, in the same way as above, the sum of
the topological indexes of $\phi ^{m_2}$ on points belonging to $K_2$ is
zero (because $\phi ^{m_2}$ is also homotopic to a Dehn twist). Again, in
the same way as above, for some adequate sufficiently large integer $m_3>0,$
the local dynamics at points in $K_2$ implies that 
$$
\stackunder{z\in K_2}{\sum }ind(\phi ^{m_3.m_2},z)<0. 
$$

So, if we project the fixed points of $\widetilde{\phi }^{m_3.m_2}$ to the
torus, we get a finite set $K_3\supset K_2,$ not equal to $K_2$ by the same
reason as above. In this way, we get a strictly increasing sequence of
finite sets $K_i.$ So at some $i_0,$ the cardinality of $K_{i_0}$ is larger
than the number of singularities of the foliations ${\cal F}^u,{\cal F}^s.$
This implies that some $\widetilde{\phi }$-periodic point $\widetilde{Q}$
does not fall into a singularity of the foliations ${\cal F}^u,{\cal F}^s$
when projected to the torus. So $Q\stackrel{def.}{=}p(\widetilde{Q})$ is
hyperbolic and its stable and unstable manifolds have a transverse
intersection in the torus. This follows from the fact that $\phi $ is
pseudo-Anosov relative to some finite set. More precisely, for any two $\phi 
$-periodic points (which may coincide), the stable manifold of one of them
intersects the unstable manifold of the other, see \cite{FLP}. 
$\Box $

\vskip 0.2truecm

Now we continue the proof of lemma \ref{hiperdehn}. Note that $\phi $ is
transitive, preserves area and by proposition \ref{ptohiperpA}, $\phi $ has
hyperbolic periodic points with transversal homoclinic intersections in the
torus for all rotation numbers in $]\rho _2,\rho _1[.$ 
The following proposition can be extracted from the proof of lemma \ref
{interhet}:

\begin{proposition}
\label{umconseta}: Assume that $\varphi :{\rm T^2\rightarrow T^2}$ is a
transitive area preserving homeomorphism of the torus, either homotopic to
the identity or to a Dehn twist, and either $(0,0)$ or $0$ belongs to the
interior of the rotation set of a (fixed) lift of $\varphi $ to the plane,
denoted $\widetilde{\varphi }.$ Assume also that $\varphi $ has hyperbolic
periodic points with transversal homoclinic intersections 
for all rotation vectors or numbers in a small neighborhood of the origin.
Then, $\widetilde{\varphi }$ has a hyperbolic periodic point denoted $
\widetilde{Q},$ such that $W^u(\widetilde{Q})$ has a topologically
transverse intersection with $W^s(\widetilde{Q}+(a,b))$ for some pair of
coprime integers $a,b$ with $(a,b)\neq (0,0).$
\end{proposition}

{\bf Remark: }When we say that $\widetilde{\varphi }$ has some hyperbolic
periodic point, we mean that it is differentiable at the whole orbit of this
point and satisfies the usual definition of hyperbolicity on the orbit.

In the same way, the following statement can be extracted from the proofs of
theorems \ref{densevarident} and \ref{densevardehn}:

\begin{proposition}
\label{doisconseta}: Under the same hypothesis as in proposition \ref
{umconseta}, the following holds: $\widetilde{\varphi }$ has some periodic
point denoted $\widetilde{Q},$ which is hyperbolic and $W^u(\widetilde{Q})$
has a topologically transverse intersection with $W^s(\widetilde{Q}+(a,b))$
for all pairs of integers $(a,b).$\ 
\end{proposition}

So $\widetilde{\phi }$ has a hyperbolic periodic point $\widetilde{Q}$ such
that $W^u(\widetilde{Q})$ has a transverse intersection with $W^s(\widetilde{%
Q}).$ The intersection is really transverse, not only topologically
transverse because $\phi $ is pseudo-Anosov rel to some finite $\phi $%
-invariant set. Thus all intersections between the unstable and stable
foliations ${\cal F}^u,{\cal F}^s$ (at regular points, which is the case)
are transverse, see \cite{FLP}.

And this gives a compact $\widetilde{\phi }$-invariant set $K_\phi \subset 
{\rm I}\negthinspace {\rm R^2}$ which carries topological entropy, that is $%
h_{top}(\widetilde{\phi }\mid _{K_\phi })>0$ ($K_\phi $ is simply an
invariant set of the horseshoe $\widetilde{\phi }$ has at the point $
\widetilde{Q}$). The only place in the proofs of lemma \ref{interhet} and
theorems \ref{densevarident} and \ref{densevardehn} where $C^{1+\epsilon }$
differentiability is necessary is in order to guarantee the existence of
hyperbolic periodic points with transversal homoclinic intersections in the
torus for rotation numbers in the interior of the vertical rotation set. As
we already said, here this is achieved by proposition \ref{ptohiperpA}.

To conclude, we need a similar statement for $\widetilde{f},$ namely we have
to show that there exists a compact $\widetilde{f}$-invariant set $%
K_f\subset {\rm I}\negthinspace {\rm R^2}$ such that $h_{top}(\widetilde{f}%
\mid _{K_f})>0.$ With this purpose, we use theorem \ref{handel} and the
remark after it. And finally we apply theorem \ref{katokok}. This part of
the proof is exactly as in end of the proof of lemma \ref{hiperident}. 
$\Box $

\vskip 0.2truecm

\begin{lemma}
\label{nodisctrans}: Suppose $f$ belongs to $Diff_0^{1+\epsilon }({\rm T^2})$
and $int(\rho (\widetilde{f}))$ is not empty or $f$ belongs to $%
Diff_k^{1+\epsilon }({\rm T^2})$ and $int(\rho _V(\widetilde{f}))$ is not
empty . If $f$ is transitive, then $f$ can not have a periodic open disk.
\end{lemma}

{\it Proof:}

The proof in both cases is analogous, so let us do it in the homotopic to
the identity case. By contradiction, suppose that for some open disk $%
D\subset {\rm T^2}$ there exists $n>0$ such that $f^n(D)=D.$ Then, as $f$ is
transitive, 
$$
Orb(D)=D\cup f(D)\cup ...\cup f^{n-1}(D) 
$$
is a $f$-invariant open set, dense in the torus. Moreover, as $f^n(D)=D,$ we
get that there exists an integer vector $(k_1,k_2)$ such that for any
connected component $\widetilde{D}$ of $p^{-1}(D),$ we have:%
$$
\widetilde{f}^n(\widetilde{f}^i(\widetilde{D}))=\widetilde{f}^i(\widetilde{D}%
)+(k_1,k_2),\text{ for all }0\leq i\leq n-1 
$$

As $int(\rho (\widetilde{f}))\neq \emptyset ,$ choose a rational vector $%
\left( \frac pq,\frac rq\right) \in int(\rho (\widetilde{f}))\backslash
\left( \frac{k_1}n,\frac{k_2}n\right) .$ From lemma \ref{hiperident}, $f$
has a hyperbolic $m.q$-periodic point $z$ (for some integer $m>0)$ which
realizes this rotation vector, with a transversal homoclinic intersection.

Now we claim that $(W^u(z)\cup z\cup W^s(z))\cap Orb(D)=\emptyset .$
Clearly, $z\notin Orb(D),$ because $\rho (\widetilde{z})=\left( \frac
pq,\frac rq\right) \neq \left( \frac{k_1}n,\frac{k_2}n\right) .$ Suppose by
contradiction, that a branch $\Gamma $ of $W^u(z)$ or of $W^s(z)$ intersects
some $f^i(D)$ (without loss of generality, we can suppose that $i=0$). Let $%
w\in \Gamma \cap D.$\ As $f^{n.m.q}(\Gamma )=\Gamma ,$ first suppose that $%
\Gamma $ is contained in $D.$ In this case, given $\widetilde{D}$ a
connected component of $p^{-1}(D),$ there exists a connected component $
\widetilde{\Gamma }\in p^{-1}(\Gamma )$ such that $\widetilde{\Gamma }%
\subset \widetilde{D}.$ As $\widetilde{f}^{n.m.q}(\widetilde{\Gamma })= 
\widetilde{\Gamma }+(n.m.p,n.m.r),$ $\widetilde{f}^{n.m.q}(\widetilde{D})= 
\widetilde{D}+(m.q.k_1,m.q.k_2)$ and $(n.m.p,n.m.r)\neq (m.q.k_1,m.q.k_2),$
we get a contradiction. So, there are points $w^{\prime },w^{\prime \prime
}\in \Gamma \cap \partial D$ such that $w$ belongs to the arc in $\Gamma $
between $w^{\prime }$ and $w^{\prime \prime }$ and apart from its end
points, this arc, denoted $\gamma $ is contained in $D$. Clearly $\gamma $
divides $D$ into two open disks, $D_1$ and $D_2.$ So, as $f^{n.m.q}\times
f^{n.m.q}:D\times D\rightarrow D\times D$ preserves volume, we get from the
Poincar\'e recurrence theorem that there exists an integer $N>0$ such that $%
f^{N.n.m.q}\times f^{N.n.m.q}(D_1\times D_2)$ intersects $D_1\times D_2.$
But this is a contradiction with the preservation of area of $f^{n.m.q}$ and
the fact that $\gamma $ is an arc contained in $\Gamma ,$ because $\Gamma $
is a branch of $W^u(z)$ or of $W^s(z).$ So $\gamma $ can not be $f$-periodic
neither can some $f$-iterate of $\gamma $ have a topologically transverse
intersection with $\gamma $$.$ This argument appears in the proof of lemma
6.1 of \cite{frankslecal}.

So $\Gamma $ does not intersect $Orb(D),$ that is, $(W^u(z)\cup z\cup
W^s(z))\subset \partial Orb(D)$ because $Orb(D)$ is dense in the torus. As $%
Orb(D)$ is an open $f$-invariant set, if $K$ is a connected component of $%
Orb(D),$ as $f$ is transitive, the first natural $n_K$ such that $%
f^{n_K}(K)=K$ (the existence of $n_K$ follows from the fact that f is area
preserving) satisfies the following 
$$
Orb(D)=K\cup f(K)\cup ...\cup f^{n_K-1}(K)\text{ and the union is disjoint. }
$$
Thus, $Orb(D)$ has exactly $n_K$ connected components. But as $(W^u(z)\cup
z\cup W^s(z))\subset \partial Orb(D)$ and $W^u(z)$ has a transversal
intersection with $W^s(z),$ we get that $Orb(D)$ has infinitely many
connected components, a contradiction. $\Box $

\vskip 0.2truecm

\begin{lemma}
\label{interhet}: Suppose $f$ belongs to $Diff_0^{1+\epsilon }({\rm T^2})$
and $(0,0)\in int(\rho (\widetilde{f}))$ or $f$ belongs to $%
Diff_k^{1+\epsilon }({\rm T^2})$ and $0\in int(\rho _V(\widetilde{f})).$ If $%
f$ is transitive and $\widetilde{Q}\in {\rm I}\negthinspace 
{\rm R^2}$ is a hyperbolic periodic point for $\widetilde{f}$ with
transversal homoclinic intersections when projected to the torus, then for
some pair of integers $(a,b)\neq (0,0),$ $a,b$ coprimes, $W^u(\widetilde{Q})%
\pitchfork W^s(\widetilde{Q}+(a,b)).$
\end{lemma}

{\it Proof:}

Again, the proof in both situations is analogous, so let us present the
arguments in the homotopic to the identity case. Let $\widetilde{Q}\in {\rm I%
}\negthinspace 
{\rm R^2}$ be a hyperbolic periodic point for $\widetilde{f}$ with a
transversal homoclinic intersection, which exists by lemma \ref{hiperident}.
First note that there exists arbitrarily small topological rectangles $%
D_Q\subset {\rm T^2}$ such that $Q=p(\widetilde{Q})$ is a vertex of $D_Q$
and the sides of $D_Q,$ denoted $\alpha _Q,\beta _Q,\gamma _Q$ and $\delta
_Q $ are contained in $W^s(Q),W^u(Q),W^s(Q)$ and $W^u(Q)$ respectively. Let
us choose another hyperbolic periodic point $P\in {\rm T^2} $ such that its
rotation vector is not $(0,0),$ which also has a transversal homoclinic
intersection (the existence of this point $P$ also follows from lemma \ref
{hiperident}). Let $n$ be a natural number such that $f^n(Q)=Q$ and $%
f^n(P)=P $ and $Df^n$ has positive eigenvalues at both points. Clearly the
orbit of $Q $ is disjoint from the orbit of $P.$ So we can also choose
another arbitrarily small topological rectangle $D_P\subset {\rm T^2}$ such
that $P$ is a vertex of $D_P$ and the sides of $D_P,$ denoted $\alpha
_P,\beta _P,\gamma _P$ and $\delta _P$ are contained in $%
W^s(P),W^u(P),W^s(P) $ and $W^u(P)$ respectively, in such a way that 
$$
\left( D_Q\cup f^{-1}(D_Q)\cup ...\cup f^{-n+1}(D_Q)\right) \cap \left(
D_P\cup f^{-1}(D_P)\cup ...\cup f^{-n+1}(D_P)\right) =\emptyset 
$$
and for all $0\leq i\leq n-1,$ each set $f^{-i}(D_Q),$ $f^{-i}(D_P)$ is a
small rectangle in ${\rm T^2.}$

Also note that there exists an integer $m_0>0$ such that for all $0\leq
i\leq n-1,$ if $m\geq m_0,$ then 
$$
\begin{array}{c}
\partial (f^{n.m}(f^{-i}(D_Q)))\subset W^u(f^{-i}(Q))\cup f^{-i}(\alpha _Q)
\\ 
\text{and} \\ \partial (f^{n.m}(f^{-i}(D_P)))\subset W^u(f^{-i}(P))\cup
f^{-i}(\alpha _P), 
\end{array}
$$
because%
$$
f^{n.m}(f^{-i}(\alpha _Q))\subset f^{-i}(\alpha _Q)\text{ and }%
f^{n.m}(f^{-i}(\alpha _P))\subset f^{-i}(\alpha _P)\text{ for all integers }%
m>0\ 
$$
and for $m\geq m_0,$ $f^{n.m}(f^{-i}(\gamma _Q))\subset f^{-i}(\alpha _Q)$
and $f^{n.m}(f^{-i}(\gamma _P))\subset f^{-i}(\alpha _P).$

\vskip 0.1truecm

As $f$ is transitive, for all $0\leq i\leq n-1,$ there exist integers $%
l_P(i),l_Q(i)\geq m_0.n$ such that 
$$
f^{l_Q(i)}(f^{-i}(D_Q))\cap D_P\neq \emptyset \text{ and }%
f^{l_P(i)}(f^{-i}(D_P))\cap D_Q\neq \emptyset . 
$$

So for any $0\leq i\leq n-1,$ there exists integers $m_P(i),m_Q(i)\geq m_0$
and other integers $0\leq r_P(i),r_Q(i)\leq n-1$ such that 
$$
f^{m_Q(i).n}(f^{-i}(D_Q))\cap f^{-r_P(i)}(D_P)\neq \emptyset \text{ and }%
f^{m_P(i).n}(f^{-i}(D_P))\cap f^{-r_Q(i)}(D_Q)\neq \emptyset . 
$$

This means that for all $0\leq i\leq n-1,$ $W^u(f^{-i}(Q))\pitchfork
W^s(f^{-r_P(i)}(P))$ and $W^u(f^{-i}(P))\pitchfork W^s(f^{-r_Q(i)}(Q)).$
Then a simple combinatorial argument implies that there exists $0\leq
i,j\leq n-1$ such that $W^u(f^{-i}(Q))\pitchfork W^s(f^{-j}(P))$ and $%
W^u(f^{-j}(P))\pitchfork
W^s(f^{-i}(Q)).$

So maybe after renaming the points in the orbit of $P,$ we can assume that $%
W^u(Q)\pitchfork W^s(P)$ and $W^u(P)\pitchfork W^s(Q).$ If we go to the
plane, we get that, fixed some $\widetilde{Q}_0\in p^{-1}(Q),$ there exists $
\widetilde{P}_0\in p^{-1}(P)$ such that $W^u(\widetilde{Q}_0)\pitchfork W^s( 
\widetilde{P}_0).$ The choice of $Q$ and $P$ implies that $\widetilde{f}^n( 
\widetilde{Q}_0)=\widetilde{Q}_0$ and there exists some pair of integers $%
(a_1,b_1)\neq (0,0)$ such that $\widetilde{f}^n(\widetilde{P}_0)=\widetilde{P%
}_0+(a_1,b_1).$ So $W^u(\widetilde{Q}_0)\pitchfork \left( W^s(\widetilde{P}%
_0+m.(a_1,b_1))\right) $ for all integers $m>0$ (clearly $W^s(\widetilde{P}%
_0+m.(a_1,b_1))=W^s(\widetilde{P}_0)+m.(a_1,b_1)$). From the fact that $%
W^u(P)\pitchfork W^s(Q)$ we get that for any $\widetilde{P}\in p^{-1}(P),$
there exists a certain $\widetilde{Q}=function(\widetilde{P}),$ such that $%
W^u(\widetilde{P})\pitchfork W^s(\widetilde{Q})$ and $\left\| \widetilde{P}- 
\widetilde{Q}\right\| <Const,$ which does not depend on the choice of $
\widetilde{P}\in p^{-1}(P).$ From the topological transversality, we get
that a piece of $W^u(Q)$ gets arbitrarily close in the Haudorff topology to
any given compact part of $W^u(P)$ (see definition 9), so if we choose an
integer $m>0$ sufficiently large, we get that a piece of $W^u(\widetilde{Q}%
_0)$ is sufficiently close to part of $W^u(\widetilde{P}_0+m.(a_1,b_1)),$
something that forces $W^u(\widetilde{Q}_0)$ to have a topological
transverse intersection with $W^s(\widetilde{Q}^{\prime })$ for some $
\widetilde{Q}^{\prime }\in p^{-1}(Q)$ such that 
$$
\left\| (\widetilde{P}_0+m.(a_1,b_1))-\widetilde{Q}^{\prime }\right\|
<Const. 
$$
Thus if $m>0$ is sufficiently large, $\widetilde{Q}^{\prime }\neq \widetilde{%
Q}_0.$ To conclude the proof we just have to note that $W^u(\widetilde{Q}%
_0)\cup \widetilde{Q}_0\cup W^s(\widetilde{Q}_0)$ is an arc connected subset
of ${\rm I}\negthinspace 
{\rm R^2,}$ which transversely intersects $(W^u(\widetilde{Q}_0)\cup 
\widetilde{Q}_0\cup W^s(\widetilde{Q}_0))+(c^{\prime },d^{\prime }),$ for
some pair of integers $(c^{\prime },d^{\prime })\neq (0,0).$ If $(c^{\prime
},d^{\prime })$ equals some integer $j$ times $(c,d),$ with $c$ and $d$
coprimes, we need the following lemma, which is a consequence of the so
called Brouwer lemma on translation arcs, see corollary 3.3 of \cite{brown}.

\begin{lemma}
\label{brouwer}:\ Let $K$ be an arc connected subset of ${\rm I}%
\negthinspace 
{\rm R^2}$ and let $\widetilde{h}:{\rm I}\negthinspace 
{\rm R^2\rightarrow I}\negthinspace 
{\rm R^2}$ be an orientation preserving homeomorphism of the plane with no
fixed points. If $\widetilde{h}(K)\cap K=\emptyset ,$ then $\widetilde{h}%
^n(K)\cap K=\emptyset $ for all integers $n\neq 0.$
\end{lemma}

The above lemma applied to $K=W^u(\widetilde{Q}_0)\cup \widetilde{Q}_0\cup
W^s(\widetilde{Q}_0)$ and $\widetilde{h}=Id+(c,d)$ implies that $W^u( 
\widetilde{Q}_0)\cup \widetilde{Q}_0\cup W^s(\widetilde{Q}_0)$ must
intersect $(W^u(\widetilde{Q}_0)\cup \widetilde{Q}_0\cup W^s(\widetilde{Q}%
_0))+(c,d),$ also in a topologically transverse way, and so the proof of
lemma \ref{interhet} is over. $\Box $

\vskip 0.2truecm

{\bf Remarks:}

\begin{enumerate}
\item  note that the above argument also implies that there exists a pair of
integers $(p,q)\neq (0,0),$ $p$ and $q$ coprimes, such that for any $
\widetilde{P}\in p^{-1}(P),$ $W^u(\widetilde{P})\pitchfork W^s(\widetilde{P}%
+(p,q)).$ As the only hypothesis on $\widetilde{P}$ was that its rotation
vector is not $(0,0),$ we get that this is a general property, true for all
rational vectors in the interior of $\rho (\widetilde{f});$

\item  As $W^u(\widetilde{Q})\pitchfork W^s(\widetilde{Q}+(a,b)),$ it is
easy to prove that $\overline{W^u(\widetilde{Q})}\supset W^u(\widetilde{Q}%
)+(a,b);$
\end{enumerate}

\begin{lemma}
\label{vardensat2}: In case $f\in Diff_0^{1+\epsilon }({\rm T^2}),$ let $P$
be a hyperbolic $m$-periodic point for $f$ with a transversal homoclinic
intersection, whose rotation vector lies in the interior of $\rho ( 
\widetilde{f})$ and in case $f\in Diff_k^{1+\epsilon }({\rm T^2}),$ let $P$
be a hyperbolic $m$-periodic point for $f$ with a transversal homoclinic
intersection, whose vertical rotation number lies in the interior of $\rho
_V(\widetilde{f})$ (in both cases, $m>0$ is an integer). If $f$ is
transitive, then $\overline{W^u(P)}=\overline{W^s(P)}={\rm T^2.}$
\end{lemma}

{\it Proof: }

The proof is analogous for $W^u(P)$ and $W^s(P),$ so let us present it for $%
W^u(P).$ As $f^m(\overline{W^u(P)})=\overline{W^u(P)},$ if we suppose by
contradiction that $\overline{W^u(P)}\neq {\rm T^2,}$ then the open set $%
\left( \overline{W^u(P)}\right) ^c$ is non-empty and invariant under $f^m.$
First, suppose that some of its connected components contains a
homotopically non trivial simple closed curve $\gamma .$ Then 
\begin{equation}
\label{interseccc}\left( \gamma \cup f^m(\gamma )\right) \cap \left(
B_\epsilon (P)\cup W^u(P)\right) =\emptyset , 
\end{equation}
for a sufficiently small $\epsilon >0.$

From remark 1 right after the proof of lemma \ref{interhet}, we get that 
\begin{equation}
\label{defpandq} 
\begin{array}{c}
\text{there exists a pair of integers }(p,q)\neq (0,0),\text{ }p\text{ and }%
q \text{ coprimes,} \\ \text{such that for any }\widetilde{P}\in
p^{-1}(P),W^u( \widetilde{P})\pitchfork (W^s(\widetilde{P})+(p,q)). 
\end{array}
\end{equation}

In the following we will analyze separately the homotopic to the identity
case and the Dehn twist case and in both we will arrive at a contradiction.

\begin{itemize}
\item  First suppose that $f\in Diff_0^{1+\epsilon }({\rm T^2}).$
\end{itemize}

Note that there exists a pair of integers $(a_1,a_2)$ such that for any $
\widetilde{P}\in p^{-1}(P),$ $\widetilde{f}^m(\widetilde{P})-(a_1,a_2)= 
\widetilde{P}$ and $(0,0)\in int(\rho (\widetilde{f}^m(\bullet
)-(a_1,a_2))). $ If we set $\widetilde{g}(\bullet )\stackrel{def.}{=} 
\widetilde{f}^m(\bullet )-(a_1,a_2),$ then $\widetilde{g}(\widetilde{P})= 
\widetilde{P}$ for any $\widetilde{P}\in p^{-1}(P)$ and so we can construct
a simple curve $\eta $ in the plane which consists of a piece of $W^u( 
\widetilde{P})$ starting at $\widetilde{P}$ until it continues through the
piece of $W^s( \widetilde{P})+(p,q)$ which is inside $B_\epsilon (\widetilde{%
P})+(p,q)$ and contains $\widetilde{P}+(p,q),$ see figure 2. Let $\theta $
be equal to $\stackunder{i\in integers}{\cup }\eta +i(p,q),$ which is a
closed connected set. Then, it is easy to see that%
$$
\theta \cap p^{-1}(\gamma )=\emptyset . 
$$
So, each connected component $\widetilde{\gamma }$ of $p^{-1}(\gamma )$ is
parallel to $\theta ,$ namely $\widetilde{\gamma }=\widetilde{\gamma }%
+(p,q). $ Moreover, the following holds:

\begin{enumerate}
\item  $\widetilde{g}^l(\theta )\subset \stackunder{i\in integers}{\cup }%
\left( \left( B_\epsilon (\widetilde{P})\cup W^u(\widetilde{P})\right)
+i.(p,q)\right) ,$ for all integers $l>0;$

\item  $\widetilde{g}^l(\theta )$ intersects $p^{-1}(\gamma )$ if $l>0$ is
sufficiently large (this happens because $\widetilde{g}^l(\theta )\cap
\theta \neq \emptyset $ and $g$ has periodic points in the torus with
rotation vector parallel to $(-q,p));$
\end{enumerate}

So $B_\epsilon (P)\cup W^u(P)$ intersects $\gamma ,$ a contradiction with
expression (\ref{interseccc}).

\begin{itemize}
\item  Now suppose $f\in Diff_k^{1+\epsilon }({\rm T^2}).$
\end{itemize}

If $\gamma $ is not homotopic to $(1,0),$ then as in this case $f^m(\gamma )$
and $\gamma $ are in different homotopic classes and $B_\epsilon (P)\cup
W^u(P)$ contains a homotopically non-trivial simple closed curve, we get
that $B_\epsilon (P)\cup W^u(P)$ either intersects $\gamma $ or $f^m(\gamma
),$ a contradiction. So we are left to consider the case when $\gamma $ is
an horizontal curve, that is, $\gamma $ is homotopic to $(1,0).$ In this
case, concerning the value of $q$ in expression (\ref{defpandq}), there are
2 possibilities:

i) $q\neq 0.$ This implies that there exists a homotopically non-trivial
closed curve contained in $B_\epsilon (P)\cup W^u(P)$ which is homotopic to $%
(m,q),$ for some integer $m.$ Note that $m$ may not be equal to $p,$ which
appears in (\ref{defpandq}). The reason for this is that we may have to
iterate the $(p,q)$ curve contained in $W^u(P)\cup P\cup W^s(P)$ many times
under $f^m$ so that its piece contained in $W^s(P)$ becomes very small, more
precisely, contained in $B_\epsilon (P)$. Nevertheless, as $q\neq 0,$ this
curve must intersect $\gamma ,$ a contradiction;

ii) $q=0.$ In this case, if we lift everything to the cylinder $S^1\times 
{\rm I}\negthinspace 
{\rm R}$, we obtain a homotopically non trivial closed curve $\widehat{%
\theta }\subset W^u(\widehat{P})\cup \widehat{P}\cup W^s(\widehat{P}),$ for
some $\widehat{P}\in S^1\times {\rm I}\negthinspace 
{\rm R}$ such that $p(\pi ^{-1}(\widehat{P}))=P\in {\rm T^2.}$ Moreover, $
\widehat{\theta }$ is obtained as the union of two simple connected arcs,
one starting at $\widehat{P},$ contained in $W^u(\widehat{P})$ and the other
contained in $W^s(\widehat{P}),$ ending at $\widehat{P}.$ As $P$ is $f$%
-periodic with period $m,$ there exists an integer $r$ such that 
$$
\widehat{f}^m(\widehat{P})=\widehat{P}+(0,r). 
$$
So if we define $\widehat{g}(\bullet )\stackrel{def.}{=}\widehat{f}%
^m(\bullet )-(0,r),$ we get that $\widehat{g}^l(\widehat{\theta })\cap 
\widehat{\theta }\neq \emptyset $ for all integers $l$ and if $l>0$ is
sufficiently large, then $\widehat{g}^l(\widehat{\theta })\subset W^u( 
\widehat{P})\cup B_\epsilon (\widehat{P})$ and it intersects $\pi
(p^{-1}(\gamma )),$ which again contradicts (\ref{interseccc}). This happens
because of two facts: a)\ if $l$ is sufficiently large, then $\widehat{g}^l
( \widehat{\theta }\cap W^s(\widehat{P}))\subset B_\epsilon ( \widehat{P});$
b) if $\widetilde{g}$ is a lift of $\widehat{g}$ to the plane, then $0\in
int(\rho _V(\widetilde{g})),$ so $\widehat{g}^l( \widehat{\theta })$ is
connected and satisfies $p_2(\widehat{g}^l( \widehat{\theta }))\stackrel{%
l\rightarrow \infty }{\rightarrow }{\rm I} \negthinspace {\rm R.}$

Thus, in both cases, all connected components of the open set $\left( 
\overline{W^u(P)}\right) ^c$ do not contain homotopically non-trivial
(considered as torus curves) closed curves.

To conclude, suppose for some connected component $U$ of $\left( \overline{%
W^u(P)}\right) ^c$, there exists a simple closed curve $\alpha $ contained
in $U$ which is not contractible in $U.$ The following argument works in
both the homotopic to the identity case and the Dehn twist case.

From what we did above, $\alpha $ is homotopically trivial as a curve in the
torus. So in the open disk bounded by $\alpha $ there are points that belong
to $\overline{W^u(P)},$ otherwise $\alpha $ would be contractible in $U.$
But this means that the whole $\overline{W^u(P)}$ is contained in this disk,
because $\overline{W^u(P)}$ is connected. And this is a contradiction with
the first thing we did, the fact that all connected components of the open
set $\left( \overline{W^u(P)}\right) ^c$ do not contain homotopically
non-trivial simple closed curves. So all simple closed curves contained in $%
U $ are contractible in $U.$ This means $U$ is an open disk, which is
periodic under $f$ (because $f$ is area-preserving). And this contradicts
lemma \ref{nodisctrans} and proves the present lemma. $\Box $

\vskip 0.2truecm

Now we are ready to prove our main results. Before presenting the proofs, we
will restate each one of them.

\subsection{Proofs of theorems 3 and 4}

{\it Proof of theorem \ref{densevarident}:}

{\bf Statement}{\it :} {\it Suppose $f$ belongs to $Diff_0^{1+\epsilon }(%
{\rm T^2}),$ it is transitive and $\left( \frac pq,\frac rq\right) \in
int(\rho (\widetilde{f})).$ Then, $\widetilde{f}^q(\bullet )-(p,r)$ has a
hyperbolic periodic point $\widetilde{Q}$ such that $\overline{W^u( 
\widetilde{Q})}=\overline{W^s(\widetilde{Q})}={\rm I}\negthinspace {\rm R^2.}
$}

\vskip 0.2truecm

We know from lemmas \ref{hiperident} and \ref{interhet} that there exists a
hyperbolic periodic point $Q\in {\rm T^2}$ such that its rotation vector is $%
\left( \frac pq,\frac rq\right) ,$ $W^u(Q)$ has a transverse intersection
with $W^s(Q)$ and for some pair of integers $(a,b)\neq (0,0),$ $a,b$
coprimes, for all $\widetilde{Q}^{\prime }\in p^{-1}(Q),$ $W^u(\widetilde{Q}%
^{\prime })\pitchfork
W^s(\widetilde{Q}^{\prime }+(a,b)).$ As $a$ and $b$ are coprimes, there are
integers $c,d$ such that 
\begin{equation}
\label{valdet}a.d-c.b=1. 
\end{equation}
In this way, for any given pair of integers $(i_{*},j_{*}),$ there exists
another pair of integers $(i,j)$ such that 
\begin{equation}
\label{sistsol}i.(a,b)+j.(c,d)=(i_{*},j_{*}). 
\end{equation}

So fixed some $\widetilde{Q}\in p^{-1}(Q),$ $W^u(\widetilde{Q}+j(c,d)+i(a,b))%
\pitchfork W^s(\widetilde{Q}+j(c,d)+(i+1)(a,b)),$ for all integers $i,j.$
Now let $\widetilde{g}(\bullet )\stackrel{def.}{=}\widetilde{f}^q(\bullet
)-(p,r)$ and suppose $m>0$ is an integer such that $\widetilde{g}^m( 
\widetilde{Q})=\widetilde{Q}.$ As in the proof of lemma \ref{vardensat2}, we
can construct a path connected set $\theta \subset {\rm I}\negthinspace {\rm %
R^2}$ such that: %
%
%
%

\begin{enumerate}
\item  $\theta =\theta +(a,b);$

\item  $\theta $ contains $\widetilde{Q}+i(a,b),$ for all integers $i;$

\item  $\theta $ contains a simple arc $\eta $ from $\widetilde{Q}$ to $
\widetilde{Q}+(a,b)$ of the following form: $\eta $ starts at $\widetilde{Q}%
, $ contains a piece $\lambda $ of a branch of $W^u(\widetilde{Q})$ until it
reaches $W^s(\widetilde{Q})+(a,b)$ and then $\eta $ contains a piece $\mu $
of a branch of $W^s(\widetilde{Q})+(a,b)$ until it reaches $\widetilde{Q}%
+(a,b),$ see figure 3. Clearly $\mu \subset W^s(\widetilde{Q})+(a,b)$ can be
chosen arbitrarily small. The set $\theta $ is obtained from $\eta $ in the
following natural way: $\theta =\stackunder{i\in integers}{\cup }\eta
+i(a,b);$

\item  $\theta $ is bounded in the $(-b,a)$ direction, that is, $\theta $ is
contained between two straight lines, both parallel to $(a,b);$ %
%
\end{enumerate}

For any given integer $l,$ let us choose two straight lines, $L_1(l)$ and $%
L_2(l)$ both parallel to $(a,b)$ such that $\theta +l(c,d)$ is contained
between them and $(L_1(l)\cup L_2(l))\cap \theta =\emptyset .$ Now we
remember that as $(0,0)\in int(\rho (\widetilde{g})),$ there exist two
periodic points for $g,$ one with rotation vector of the form $(\frac{-b}%
N,\frac aN)$ and the other with rotation vector of the form $(\frac bN,\frac{%
-a}N),$ for some sufficiently large integer $N>0,$ see \cite{franksrat}. So
as $\widetilde{g}^{m.t}(\theta )\cap \theta \neq \emptyset $ for all
integers $t,$ there exists an integer $t_{*}>0$ such that $\widetilde{g}%
^{m.t}(\theta )$ intersects both $L_1(l)$ and $L_2(l)$ for all $t\geq t^{*}.$
The way $\theta $ is constructed implies that %
%
%
there exists an integer $i_l$ such that 
\begin{equation}
\label{quasetudo}W^u(\widetilde{Q})\pitchfork W^s(\widetilde{Q}%
+l(c,d)+i(a,b)),\text{ for all }i\geq i_l. 
\end{equation}

To prove that for any given pair of integers $(i_{*},j_{*}),$ $W^u( 
\widetilde{Q})\pitchfork
W^s(\widetilde{Q}+(i_{*},j_{*})),$ from expressions (\ref{valdet}), (\ref
{sistsol}) and (\ref{quasetudo}) it is enough to show that $W^u(\widetilde{Q}%
)\pitchfork W^s(\widetilde{Q}-(a,b)).$ And this is achieved by an argument
similar to the one in the proof of lemma \ref{interhet}. To be precise, in
that proof we can choose the hyperbolic periodic point $P\in {\rm T^2}$ such
that its rotation vector (for the map $\widetilde{g}$) is $(\frac{-a}{%
n^{\prime }},\frac{-b}{n^{\prime }}),$ for some integer $n^{\prime }>0$ such
that $(\frac{-a}{n^{\prime }},\frac{-b}{n^{\prime }})\in int(\rho ( 
\widetilde{g})).$ Clearly, the rotation vector of $P$ with respect to $
\widetilde{f}$ is%
$$
\left( \frac{n^{\prime }p-a}{n^{\prime }q},\frac{n^{\prime }r-b}{n^{\prime }q%
}\right) . 
$$
So, maybe after renaming the points in the orbit of $P,$ we can assume that $%
W^u(Q)\pitchfork W^s(P)$ and $W^u(P)\pitchfork W^s(Q).$ If we go back to the
plane we get that there exists $\widetilde{P}\in p^{-1}(P)$ such that $W^u( 
\widetilde{Q})\pitchfork W^s(\widetilde{P}).$

Let $n>0$ be an integer such that $g^n(Q)=Q,$ $g^n(P)=P$ and $Dg^n$ has
positive eigenvalues at both points. The choice of $Q$ and $P$ implies that $
\widetilde{g}^n(\widetilde{Q})=\widetilde{Q}$ and there exists an integer $%
c>0$ such that $\widetilde{g}^n(\widetilde{P})=\widetilde{P}+c(-a,-b).$ So $%
W^u(\widetilde{Q})\pitchfork W^s(\widetilde{P}+t.c(-a,-b))$ for all integers 
$t>0.$ From the fact that $W^u(P)\pitchfork W^s(Q)$ we get that for any $
\widetilde{P}^{\prime }\in p^{-1}(P),$ there exists a certain $\widetilde{Q}%
^{\prime }=function(\widetilde{P}^{\prime })$ such that $W^u(\widetilde{P}%
^{\prime })\pitchfork W^s(\widetilde{Q}^{\prime })$ and $\left\| \widetilde{P%
}^{\prime }-\widetilde{Q}^{\prime }\right\| <Const.$ Moreover, $Const>0$
does not depend on the choice of $\widetilde{P}^{\prime }\in p^{-1}(P).$
From the topological transversality, we get that $\overline{W^u(Q)}$
contains $W^u(P),$ so if we choose an integer $t>0$ sufficiently large, we
get that a connected piece of $W^u(\widetilde{Q})$ is sufficiently close in
the Hausdorff topology to part of $W^u(\widetilde{P}+t.c(-a,-b)),$ something
that forces $W^u(\widetilde{Q})$ to have a topological transverse
intersection with $W^s(\widetilde{Q}^{\prime })$ for some $\widetilde{Q}%
^{\prime }\in p^{-1}(Q)$ such that $\left\| (\widetilde{P}+t.c(-a,-b))- 
\widetilde{Q}^{\prime }\right\| <Const.$ Thus $\widetilde{Q}^{\prime }= 
\widetilde{Q}+t.c(-a,-b)+e_1(a,b)+e_2(c,d),$ for some integer pair $%
(e_1,e_2) $ such that $\max \{\left| e_1\right| ,\left| e_2\right| \}$ is
bounded independently of $t>0.$ Now using expression (\ref{quasetudo}) we
get that if $t>0$ is sufficiently large, then $W^u(\widetilde{Q}^{\prime })%
\pitchfork W^s(\widetilde{Q}-i(a,b)),$ for some $i>0.$ Thus $W^u(\widetilde{Q%
})\pitchfork W^s(\widetilde{Q}-(a,b)).$

To conclude our proof, note that lemma \ref{vardensat2} implies that given
an open ball $U\subset {\rm I}\negthinspace {\rm R^2,}$ there exists $
\widetilde{Q}^{\prime }\in p^{-1}(Q)$ such that $W^u(\widetilde{Q}^{\prime
})\cap U\neq \emptyset .$ As $W^u(\widetilde{Q})\pitchfork W^s(\widetilde{Q}%
^{\prime }),$ we get that $W^u(\widetilde{Q})\cap U\neq \emptyset $ and we
are done. The proof for $W^s(\widetilde{Q})$ is analogous. $\Box $

\vskip 0.2truecm

{\it Proof of theorem \ref{densevardehn}:}

{\bf Statement}{\it :} {\it Suppose $f$ belongs to $Diff_k^{1+\epsilon }(%
{\rm T^2}),$ it is transitive, $\frac pq\in int(\rho _V(\widetilde{f}))$ and 
$s$ is any integer number. Then, $\widetilde{f}^q(\bullet )-(s,p)$ has a
hyperbolic periodic point $\widetilde{Q}$ such that $\overline{W^u( 
\widetilde{Q})}=\overline{W^s(\widetilde{Q})}={\rm I}\negthinspace {\rm R}$}$%
{\rm ^2.}$

\vskip 0.2truecm

Let $f\in Diff_{k_{Dehn}}^{1+\epsilon }({\rm T^2})$ for some integer $%
k_{Dehn}\neq 0.$ Without loss of generality, let us suppose that $%
k_{Dehn}>0. $ We know from lemma \ref{interhet} that there exists a
hyperbolic periodic point $\widetilde{Q}\in {\rm I}\negthinspace 
{\rm R^2}$ for $\widetilde{g}(\bullet )\stackrel{def.}{=}\widetilde{f}%
^q(\bullet )-(s,p),$ of period $m>0,$ such that for some pair of integers $%
(a,b)\neq (0,0),$ $a,b$ coprimes, $W^u(\widetilde{Q})\pitchfork W^s( 
\widetilde{Q}+(a,b)).$

First, let us suppose that $b>0.$ As $\widetilde{g}^m(\widetilde{Q})= 
\widetilde{Q}$ and $g$ is homotopic to a Dehn twist, we get that for any
integer $l,$ $\widetilde{g}^m(\widetilde{Q}+(0,l))=\widetilde{Q}%
+(k_{Dehn}.q.m.l,l).$ So, $W^u(\widetilde{Q})\pitchfork W^s(\widetilde{Q}%
+(a+k_{Dehn}.m.q.b,b))$ and there exists a compact path connected set $%
\gamma \subset {\rm I}\negthinspace 
{\rm R^2}$ connecting $\widetilde{Q}$ to $\widetilde{Q}+(k_{Dehn}.m.q.b,0)$
of the following form. The set $\gamma $ is contained in 
$$
W^u(\widetilde{Q})\cup W^s(\widetilde{Q}+(a+k_{Dehn}.m.q.b,b))\cup W^u( 
\widetilde{Q}+(k_{Dehn}.m.q.b,0)) 
$$
and it contains a connected arc of $W^s(\widetilde{Q}+(a+k_{Dehn}.m.q.b,b))$
that has one endpoint in $\widetilde{Q}+(a+k_{Dehn}.m.q.b,b).$

Thus $\pi (\gamma )\stackrel{def.}{=}\widehat{\gamma }\subset S^1\times {\rm %
I}\negthinspace 
{\rm R}$ contains a homotopically non trivial simple closed curve and $
\widehat{\gamma }\subset W^u(\widehat{Q})\cup W^s(\widehat{Q}+(0,b)),$ where 
$\widehat{Q}=\pi (\widetilde{Q}).$ As $\widehat{\gamma }$ is compact, for
every integer $n$ there exist numbers $M_{-}(n)<M_{+}(n)$ such that:

\begin{enumerate}
\item  $\widehat{\gamma }+(0,n)\subset S^1\times ]M_{-}(n),M_{+}(n)[;$

\item  $S^1\times \{M_{-}(n),M_{+}(n)\}$ does not intersect $\widehat{\gamma 
};$
\end{enumerate}

As $\widehat{g}$ has points with positive and negative vertical rotation
number and $\widehat{g}^{i.m}(\widehat{\gamma })\cap \widehat{\gamma }\neq
\emptyset $ for all integers $i,$ we get that for every integer $n,$ there
exists an integer $i(n)>0$ such that 
$$
\widehat{g}^{i.m}(\widehat{\gamma })\text{ intersects }S^1\times
\{M_{-}(n)\} \text{ and }S^1\times \{M_{+}(n)\}\text{ for all }i>i(n). 
$$
This means that $W^u(\widehat{Q})\pitchfork W^s(\widehat{Q}+(0,n+b))$ for
all integers $n.$ So, $W^u(\widetilde{Q})\pitchfork W^s(\widetilde{Q}%
+(c(l),l))$ for all integers $l$ and for some function $c(l).$

If we remember that $\widetilde{g}^{m.i}(\widetilde{Q}+(0,l))=\widetilde{Q}%
+(k_{Dehn}.q.m.i.l,l)$ for any integers $l$ and $i>0,$ we get that $W^u( 
\widetilde{Q})\pitchfork W^s(\widetilde{Q}+(c(l)+k_{Dehn}.q.m.i.l,l))$ for
all integers $l$ and $i>0.$ Using this it is easy to see that for some
integer constants $c_{-}<0<c_{+},$ the following intersections hold: $W^u( 
\widetilde{Q})\pitchfork W^s(\widetilde{Q}+(c_{+},0))$ and $W^u(\widetilde{Q}%
)\pitchfork W^s(\widetilde{Q}+(c_{-},0)).$ And finally, the same argument
applied in the end of the proof of lemma \ref{interhet} using lemma \ref
{brouwer} implies that $W^u(\widetilde{Q})\pitchfork W^s(\widetilde{Q}%
+(1,0)) $ and $W^u(\widetilde{Q})\pitchfork W^s(\widetilde{Q}-(1,0)).$

If $b=0,$ there exists a continuous arc $\gamma \subset {\rm I}%
\negthinspace 
{\rm R^2}$ connecting $\widetilde{Q}$ to $\widetilde{Q}+(a,0)$ contained in $%
W^u(\widetilde{Q})\cup W^s(\widetilde{Q}+(a,0)).$ Thus $\pi (\gamma )= 
\widehat{\gamma }\subset S^1\times {\rm I}\negthinspace 
{\rm R}$ contains a homotopically non trivial simple closed curve and $\pi
(\gamma )\subset W^u(\widehat{Q})\cup W^s(\widehat{Q}),$ where $\widehat{Q}%
=\pi (\widetilde{Q}).$ Arguing exactly as in the $b>0$ case we get that $%
W^u( \widehat{Q})\pitchfork W^s(\widehat{Q}+(0,n))$ for all integers $n.$
So, $W^u(\widetilde{Q})\pitchfork W^s(\widetilde{Q}+(c(n),n))$ for all
integers $n $ and for some function $c(n).$ And finally we obtain, again
exactly as in the $b>0$ case, that $W^u(\widetilde{Q})\pitchfork W^s( 
\widetilde{Q}+(1,0))$ and $W^u(\widetilde{Q})\pitchfork W^s(\widetilde{Q}%
-(1,0)).$ %
%
%
%
%
%
%
%

The case $b<0$ is analogous to the case $b>0.$

In this way, given a pair of integers $(c,d),$ we know that for some integer 
$s,$ $W^u(\widetilde{Q})\pitchfork W^s(\widetilde{Q}+(s,d)).$ But this
implies that $W^u(\widetilde{Q}+(c-s,0))\pitchfork W^s(\widetilde{Q}+(c,d)).$
As $W^u(\widetilde{Q})\pitchfork W^s(\widetilde{Q}+(c-s,0)),$ we finally get
that $W^u(\widetilde{Q})\pitchfork W^s(\widetilde{Q}+(c,d)).$ As in the
previous theorem, note that lemma \ref{vardensat2} implies that given an
open ball $U\subset {\rm I}\negthinspace {\rm R^2,}$ there exists $
\widetilde{Q}^{\prime }\in p^{-1}(Q)$ such that $W^u(\widetilde{Q}^{\prime
})\cap U\neq \emptyset .$ As $W^u(\widetilde{Q})\pitchfork W^s(\widetilde{Q}%
^{\prime }+(a,0)),$ for all integers $a,$ we get that $W^u(\widetilde{Q}%
)\cap U\neq \emptyset $ and the proof is over. The proof for $W^s(\widetilde{%
Q})$ is analogous. $\Box $

\vskip 0.2truecm

\subsection{Proof of corollary 5}

{\bf Statement}{\it :} {\it Suppose $f$ belongs to $Diff_0^{1+\epsilon }(%
{\rm T^2})$ and $(0,0)\in int(\rho (\widetilde{f}))$ or $f$ belongs to $%
Diff_k^{1+\epsilon }({\rm T^2})$ and $0\in int(\rho _V(\widetilde{f})).$ If $%
f$ is transitive, then $\widetilde{f}$ is topologically mixing.}

\vskip 0.2truecm

Remember that topologically mixing means the following: given $U,V$ be
arbitrarily small open balls in ${\rm I}\negthinspace {\rm R^2,\ }$there
exists an integer $N=N(U,V)>0$ such that if $n\geq N,$ then $\widetilde{f}%
^n(U)\cap V\neq \emptyset .$

First note that, from the previous theorems, $\widetilde{f}$ has a
hyperbolic $m$-periodic point $\widetilde{Q}$ such that $\overline{W^u( 
\widetilde{Q})}=\overline{W^s(\widetilde{Q})}={\rm I}\negthinspace {\rm R^2}$
(this is all we need to get the topological mixing). So, $W^u(\widetilde{f}%
^i(\widetilde{Q}))$ intersects $V$ and $W^s(\widetilde{f}^i(\widetilde{Q}))$
intersects $U,$ for all $0\leq i\leq m-1.$ Let $\lambda _i^s$ be a compact
connected piece of a branch of $W^s(\widetilde{f}^i(\widetilde{Q}))$ such
that $\lambda _i^s$ starts at $\widetilde{f}^i(\widetilde{Q})$ and
intersects $U,$ for all $0\leq i\leq m-1.$ Clearly, there exists $N>0,$ a
large integer, such that if $n\geq N,$ then a connected piece of $\widetilde{%
f}^{-n}(V)$ is sufficiently close in the Hausdorff topology to $\lambda
_{i(n)}^s,$ for some $0\leq i(n)\leq m-1,$ in a way that this forces $
\widetilde{f}^{-n}(V)$ to intersect $U.$ And so $V\cap \widetilde{f}%
^n(U)\neq \emptyset .$ $\Box $

\vskip 0.2truecm

\subsection{Proof of lemma 1}

%
%
%

{\bf Statement}{\it :} {\it Suppose $f$ belongs to $Diff_0^{1+\epsilon }(%
{\rm T^2})$ and $(0,0)\in int(\rho (\widetilde{f}))$ or $f$ belongs to $%
Diff_k^{1+\epsilon }({\rm T^2})$ and $0\in int(\rho _V(\widetilde{f})).$
Then there exists $\widetilde{Q}\in {\rm I}\negthinspace 
{\rm R^2,}$ which is a hyperbolic periodic point for $\widetilde{f}$ such
that for some pair of integers $(a,b)\neq (0,0),$ $a,b$ coprimes, $W^u( 
\widetilde{Q})\pitchfork W^s(\widetilde{Q}+(a,b))$} $($note that $W^s( 
\widetilde{Q}+(a,b))=W^s(\widetilde{Q})+(a,b)${\it $).$ In particular, if $%
Q=p(\widetilde{Q})$ then $W^s(Q)\cup Q\cup W^u(Q)$ contains a homotopically
non-trivial simple closed curve in ${\rm T^2.}$}

\vskip 0.2truecm

This proof is analogous in both cases, so let us suppose that $f$ belongs to 
$Diff_0^{1+\epsilon }({\rm T^2})$ and $(0,0)\in int(\rho (\widetilde{f})).$
As we did in lemma \ref{hiperident}, choose three periodic orbits $Q_1,Q_2$
and $Q_3$ as in theorem \ref{llibremack}, such that their rotation vectors
form a triangle $\Delta $ that contains $(0,0)$ in its interior. Theorem \ref
{llibremack} tells us that the isotopy class of $f$ relative to $Q_1\cup
Q_2\cup Q_3$ contains a pseudo-Anosov map $\phi :{\rm T^2\rightarrow T^2}$
rel $Q_1\cup Q_2\cup Q_3$.

As $\phi :{\rm T^2\rightarrow T^2}$ is a homeomorphism of the torus
homotopic to the identity, there exists a lift of $\phi $ to the plane,
denoted $\widetilde{\phi }:{\rm I}\negthinspace {\rm R^2\rightarrow I}%
\negthinspace {\rm R^2}$ which is equivariantly homotopic to $\widetilde{f}$
rel to $p^{-1}(Q_1\cup Q_2\cup Q_3)$ and its rotation set $\rho (\widetilde{%
\phi })\supset \Delta \supset int(\Delta )\ni (0,0).$

\begin{proposition}
\label{pAdense}: The map $\widetilde{\phi }$ is topologically mixing in $%
{\rm I}\negthinspace 
{\rm R^2,}$ which implies that for every integer $N>0,$ $\widetilde{\phi }^N$
is a transitive map of the plane.
\end{proposition}

{\it Proof:}

In the proof of lemma \ref{hiperdehn} we presented propositions \ref
{umconseta} and \ref{doisconseta}, which can be directly extracted from the
proofs of lemma \ref{interhet} and theorems \ref{densevarident} and \ref
{densevardehn}. Analogously, the next proposition can be extracted from the
proof of lemma \ref{vardensat2}:

\begin{proposition}
\label{tresconseta}: Assume that $\varphi :{\rm T^2\rightarrow T^2}$ is a
transitive area preserving homeomorphism of the torus, either homotopic to
the identity or to a Dehn twist, such that for some (fixed) lift of $\varphi 
$ to the plane, denoted $\widetilde{\varphi },$ an open neighborhood of
zero, denoted $V,$ belongs to the interior of its rotation set. Also assume
that $\varphi $ has hyperbolic periodic points with transversal homoclinic
intersections for all rotation vectors (or numbers) in $V.$ If $Q\in {\rm T^2%
}$ is a hyperbolic $\varphi $-periodic point with transversal homoclinic
intersections, whose rotation vector (or number) belongs to $V,$ then $
\overline{W^u(Q)}=\overline{W^s(Q)}={\rm T^2.}$
\end{proposition}

As $\phi $ is transitive, preserves area and for all rationals $\rho \in
int(\Delta ),$ $\phi $ has a hyperbolic (remember that in this context,
hyperbolic means four sectors) periodic point with transverse homoclinic
intersections in the torus whose rotation vector is $\rho $ (see proposition 
\ref{ptohiperpA} and remark 2 right after it), we get from propositions \ref
{umconseta}, \ref{doisconseta} and \ref{tresconseta} that there exists a
hyperbolic $\widetilde{\phi }$-periodic point $\widetilde{Q}\in {\rm I}%
\negthinspace 
{\rm R^2}$ such that $\overline{W^u(\widetilde{Q})}=\overline{W^s(\widetilde{%
Q})}={\rm I}\negthinspace 
{\rm R^2}$ and so corollary \ref{topmix} finally proves the present
proposition. $\Box $

\vskip 0.2truecm

{\bf Remark: }The above proof also works in case $f$ is homotopic to a Dehn
twist.

From theorem \ref{handel} and the remark right after it, for any fixed
integer $N>0,$ we get that 
\begin{equation}
\label{asasas}s\circ f^N\mid _W=\phi ^N\circ s\text{ and }\widetilde{s}\circ 
\widetilde{f}^N\mid _{p^{-1}(W)}=\widetilde{\phi }^N\circ \widetilde{s} 
\end{equation}
for a certain lift $\widetilde{s}:p^{-1}(W)\rightarrow {\rm I}%
\negthinspace 
{\rm R^2}$ such that $\stackunder{\widetilde{z}\in p^{-1}(W)}{\sup }\left\| 
\widetilde{s}(\widetilde{z})-\widetilde{z}\right\| <\infty .$ As $\widetilde{%
\phi }^N$ is a transitive map of the plane (by proposition \ref{pAdense}),
there exists a point $\widetilde{z}^{*}\in p^{-1}(W)$ such that $\widetilde{s%
}(\widetilde{z}^{*})$ has a dense orbit under iterates of $\widetilde{\phi }%
^N.$ This is equivalent to saying that the $\omega $-limit set of $
\widetilde{s}(\widetilde{z}^{*})$ under $\widetilde{\phi }^N$ is the whole
plane. Expression (\ref{asasas}) then implies that 
$$
\widetilde{s}(\omega \text{-limit set of }\widetilde{z}^{*}\text{ under } 
\widetilde{f}^N)=(\omega \text{-limit set of }\widetilde{s}(\widetilde{z}%
^{*})\text{ under }\widetilde{\phi }^N)={\rm I}\negthinspace 
{\rm R^2.} 
$$

So as $\widetilde{s}\mid _{p^{-1}(W^{\prime })}$ is one to one, $\overline{%
p^{-1}(W^{\prime })}=p^{-1}(W)$ and $\widetilde{s}(p^{-1}(W^{\prime
}))=p^{-1}(\{periodic$ $points$ $of$ $\phi \}),$ which is a dense subset of
the plane (again, see the remark right after theorem \ref{handel}), we get
that the $\omega $-limit set of $\widetilde{z}^{*}$ under $\widetilde{f}^N$
is the whole $p^{-1}(W).$ In other words, we have proved the following
proposition:

\begin{proposition}
\label{consertnoite}: For every integer $N>0,$ there exists a point $
\widetilde{z}^{*}\in p^{-1}(W)$ such that its orbit under $\widetilde{f}^N$
is dense in $p^{-1}(W).$
\end{proposition}

\vskip 0.2truecm

Clearly, $p^{-1}(W)\supset K_f,$ where $K_f\subset {\rm I}\negthinspace {\rm %
R^2}$ is a $\widetilde{f}$-invariant compact set such that $h_{top}( 
\widetilde{f}\mid _{K_f})>0,$ see the proof of lemma \ref{hiperident}. Let
us fix some $\epsilon >0.$ As we already explained in subsection 2.2, there
exists a hyperbolic ergodic $\widetilde{f}$-invariant Borel probability
measure $\mu ^{*}$ with positive entropy, whose support is contained in $%
K_f. $ Now we describe two results from \cite{andre} on Pesin theory.

For every $\delta >0,$ there exists a compact subset $\Lambda _\delta
\subset K_f$ with $\mu ^{*}(\Lambda _\delta )>1-\delta ,$ which is the so
called $\delta $-Pesin set. Each of its points satisfies some local
hyperbolicity assumptions, see for instance section 3 of the supplement by
Katok-Mendoza in \cite{livrokatok}. Given any point $\widetilde{w}\in
\Lambda _\delta ,$ there exists a sufficiently small compact neighborhood $V$
of $\widetilde{w}$ such that the local stable manifolds of points $
\widetilde{w}^{\prime }\in V\cap \Lambda _\delta $ are either disjoint or
equal and they depend continuously on the point $\widetilde{w}^{\prime }\in
V\cap \Lambda _\delta .$ The same holds for the unstable manifolds. There is
also a continuous product structure: given two points in $V\cap \Lambda
_\delta $ the intersection between the stable manifold of one of them with
the unstable manifold of the other is transverse and thus consists of
exactly one point, which varies continuously with the first two points and
may not be in $\Lambda _\delta .$ There are also some new technical
definitions in \cite{andre}, called inaccessibility and accessibility of a
point in $V\cap \Lambda _\delta \cap Recurrent^{\pm },$where $Recurrent^{\pm
}$ is the subset of all points in $K_f$ which are both forward and backward
recurrent, see definition 6 of \cite{andre}. With this said, we can state
lemmas 5 and 6 of \cite{andre}:

\vskip 0.2truecm

{\bf Lemma 5 of \cite{andre}}{\it : Let $\widetilde{w}^{\prime }\in V\cap
\Lambda _\delta \cap Recurrent^{\pm }$ be an inaccessible point. Then there
exist rectangles enclosing $\widetilde{w}^{\prime }$ having sides along the
invariant manifolds of two hyperbolic saddles in $V$ and having arbitrarily
small diameter.}

\vskip 0.2truecm

{\bf Lemma 6 of \cite{andre}}{\it : The set of accessible points in $V\cap
\Lambda _\delta \cap Recurrent^{\pm }$ has zero $\mu ^{*}$ measure.}

\vskip 0.2truecm

As $\mu ^{*}(Recurrent^{\pm })=1$ and we can start with some point $
\widetilde{w}\in \Lambda _\delta \cap supp(\mu ^{*}),$ the following holds: $%
\mu ^{*}(V\cap \Lambda _\delta \cap Recurrent^{\pm })>0,$ so lemmas 5 and 6
of \cite{andre} then say that there exists a point $\widetilde{z}\in K_f$
such that arbitrarily small rectangles enclosing $\widetilde{z}$ can be
obtained, having sides along the invariant manifolds of two hyperbolic $
\widetilde{f}$-periodic points $\widetilde{Q}_{V1},\widetilde{Q}_{V2},$
whose orbits are contained in $V_\epsilon (K_f),$ see figure 4 and theorem 
\ref{katokok}. Clearly, if we want smaller rectangles, the points $
\widetilde{Q}_{V1},\widetilde{Q}_{V2}$ change, getting closer and closer to $
\widetilde{z}.$ Denote these rectangles by $\widetilde{Ret}(\widetilde{z}, 
\widetilde{Q}_{V1},\widetilde{Q}_{V2}).$ Standard arguments in Pesin theory
imply that the invariant manifolds of these periodic points have transverse
intersections. So, 
\begin{equation}
\label{varum}\overline{W^s(\widetilde{Q}_{V1})}=\overline{W^s(\widetilde{Q}%
_{V2})}\text{ and }\overline{W^u(\widetilde{Q}_{V1})}=\overline{W^u( 
\widetilde{Q}_{V2})}. 
\end{equation}

Suppose $\widetilde{Ret}(\widetilde{z},\widetilde{Q}_{V1},\widetilde{Q}%
_{V2}) $ is sufficiently small in a way that for all pairs of integers $%
(a,b)\neq (0,0),$ 
\begin{equation}
\label{vardois}\widetilde{Ret}(\widetilde{z},\widetilde{Q}_{V1},\widetilde{Q}%
_{V2})\cap (\widetilde{Ret}(\widetilde{z},\widetilde{Q}_{V1},\widetilde{Q}%
_{V2})+(a,b))=\emptyset . 
\end{equation}

Let us fix an integer $N>0$ which is a common period for $\widetilde{Q}_{V1}$
and $\widetilde{Q}_{V2}.$ For every pair of integers $(a,b),$ as the
interiors of both rectangles, $\widetilde{Ret}(\widetilde{z},\widetilde{Q}%
_{V1},\widetilde{Q}_{V2})$ and $\widetilde{Ret}(\widetilde{z},\widetilde{Q}%
_{V1},\widetilde{Q}_{V2})+(a,b),$ intersect $p^{-1}(W),$ we get from
proposition \ref{consertnoite} that there exists an integer $k(a,b)>0$ such
that 
\begin{equation}
\label{varinters}\widetilde{f}^{k(a,b).N}(\widetilde{Ret}(\widetilde{z}, 
\widetilde{Q}_{V1},\widetilde{Q}_{V2}))\cap (\widetilde{Ret}(\widetilde{z}, 
\widetilde{Q}_{V1},\widetilde{Q}_{V2})+(a,b))\neq \emptyset . 
\end{equation}

So from expressions (\ref{varum}), (\ref{vardois}) and (\ref{varinters}) we
get, that $W^u(\widetilde{Q}_{V1})\pitchfork
W^s(\widetilde{Q}_{V1}+(a,b)).$ $\Box $

\vskip 0.2truecm

\subsection{Proof of theorem 6}

{\bf Statement}{\it :} {\it Suppose $f$ belongs to $Diff_0^{1+\epsilon }(%
{\rm T^2})$ and $(0,0)\in int(\rho (\widetilde{f}))$ or $f$ belongs to $%
Diff_k^{1+\epsilon }({\rm T^2})$ and $0\in int(\rho _V(\widetilde{f})).$
Then, $\widetilde{f}$ has a hyperbolic periodic point $\widetilde{Q}$ such
that for any pair of integers $(a,b),$ $W^u(\widetilde{Q})\pitchfork W^s( 
\widetilde{Q}+(a,b)),$ so $\overline{W^u(\widetilde{Q})}=\overline{W^u( 
\widetilde{Q})}+(a,b)$ and $\overline{W^s(\widetilde{Q})}=\overline{W^s( 
\widetilde{Q})}+(a,b).$ Moreover, $\overline{W^u(\widetilde{Q})}=\overline{%
W^s(\widetilde{Q})},$ this set is $\widetilde{f}$-invariant ($\widetilde{f}( 
\overline{W^u(\widetilde{Q})})=\overline{W^u(\widetilde{Q})}$) and all
connected components of the complement of $\overline{W^u(\widetilde{Q})}$
are open disks with diameter uniformly bounded from above, $f$-periodic when
projected to the torus.}

\vskip 0.2truecm

This proof is analogous in both cases, so let us suppose that $f$ belongs to 
$Diff_0^{1+\epsilon }({\rm T^2})$ and $(0,0)\in int(\rho (\widetilde{f})).$

The first part of the theorem follows from lemma \ref{rothorse} which says
that there exists a hyperbolic $\widetilde{f}$-periodic point $\widetilde{Q}%
_0$ of period $m>0,$ such that for all pairs of integers $(a,b),$ $W^u( 
\widetilde{Q}_0)\pitchfork W^s(\widetilde{Q}_0+(a,b)).$

So $\overline{W^u(\widetilde{Q}_0)}$ is invariant under integer translations
and the same holds for $\overline{W^s(\widetilde{Q}_0)}.$ In particular,
this means that if $Q_0=p(\widetilde{Q}_0),$ then 
\begin{equation}
\label{ulttimma}\overline{W^u(\widetilde{Q}_0)}=p^{-1}(\overline{W^u(Q_0)}) 
\text{ and }\overline{W^s(\widetilde{Q}_0)}=p^{-1}(\overline{W^s(Q_0)}). 
\end{equation}

Now let us prove that $\overline{W^u(\widetilde{Q}_0)}=\overline{W^s( 
\widetilde{Q}_0)}$ and each connected component of $\left( \overline{W^u( 
\widetilde{Q}_0)}\right) ^c$ is a disk with diameter uniformly bounded from
above. We need another lemma.

\begin{lemma}
\label{discgeral}: Suppose $f$ belongs to $Diff_0^{1+\epsilon }({\rm T^2})$
and $int(\rho (\widetilde{f}))$ is not empty or $f$ belongs to $%
Diff_k^{1+\epsilon }({\rm T^2})$ and $int(\rho _V(\widetilde{f}))$ is not
empty. Then, $f$ can not have a periodic unbounded open disk. Moreover, all
periodic open disks have diameter uniformly bounded from above.
\end{lemma}

{\it Proof:}

This proof will be based on the proof of lemma \ref{nodisctrans} and it will
be presented in case $f$ belongs to $Diff_0^{1+\epsilon }({\rm T^2})$ and $%
int(\rho (\widetilde{f}))$ is not empty because a similar argument works in
case $f$ is homotopic to a Dehn twist.

Suppose that for some open disk $D\subset {\rm T^2}$ there exists $n>0$ such
that $f^n(D)=D.$ Then there exists a integer vector $(k_1,k_2)$ such that
for any connected component $\widetilde{D}$ of $p^{-1}(D),$ we have:%
$$
\widetilde{f}^n(\widetilde{f}^i(\widetilde{D}))=\widetilde{f}^i(\widetilde{D}%
)+(k_1,k_2),\text{ for all }0\leq i\leq n-1 
$$

Choose a rational vector $\left( \frac pq,\frac rq\right) \in int(\rho ( 
\widetilde{f}))\backslash \left( \frac{k_1}n,\frac{k_2}n\right) .$ From what
we already proved in theorem \ref{casogeral}, $\widetilde{f}^q(\bullet
)-(p,r)$ has a hyperbolic $m$-periodic point $\widetilde{Q}$ (for some
integer $m>0)$ such that $W^u(\widetilde{Q})\pitchfork W^s(\widetilde{Q}%
+(1,0))$ and $W^u(\widetilde{Q})\pitchfork W^s(\widetilde{Q}+(0,1)).$ Now
let us consider the following curves in the plane:

\begin{enumerate}
\item  let $\alpha _H$ be a simple curve connecting $\widetilde{Q}$ to $
\widetilde{Q}+(1,0),$ contained in $W^u(\widetilde{Q})\cup W^s(\widetilde{Q}%
+(1,0))$ in such a way that its intersection with $W^u(\widetilde{Q})$ and
also with $W^s(\widetilde{Q}+(1,0))$ is connected. Let $\theta _H=%
\stackunder{i\in integers}{\cup }\alpha _H+(i,0);$

\item  let $\alpha _V$ be a simple curve connecting $\widetilde{Q}$ to $
\widetilde{Q}+(0,1),$ contained in $W^u(\widetilde{Q})\cup W^s(\widetilde{Q}%
+(0,1))$ in such a way that its intersection with $W^u(\widetilde{Q})$ and
also with $W^s(\widetilde{Q}+(0,1))$ is connected. Let $\theta _V=%
\stackunder{i\in integers}{\cup }\alpha _V+(0,i);$
\end{enumerate}

Finally, let 
$$
K=\left( \stackunder{i\in integers}{\cup }\theta _H+(0,i)\right) \cup \left( 
\stackunder{i\in integers}{\cup }\theta _V+(i,0)\right) . 
$$

If $\widetilde{D}$ is a connected component of $p^{-1}(D),$ then there
exists a constant $Max>0$ such that if $diam(\widetilde{D})>Max,$ then $%
p(K)\cap D\neq \emptyset .$ So suppose $diam(\widetilde{D})>Max.$ As $\left(
\frac pq,\frac rq\right) \neq \left( \frac{k_1}n,\frac{k_2}n\right) ,$ $Q=p( 
\widetilde{Q})\notin D.$ So $W^u(Q)\cup W^s(Q)$ intersects $D.$ Suppose some
point $w\in p(K)\cap W^u(Q)$ belongs to $D.$\ Let $\Gamma $ be the branch of 
$W^u(Q)$ that contains $w.$ As we did in the proof of lemma \ref{nodisctrans}%
, there are points $w^{\prime },w^{\prime \prime }\in \Gamma \cap \partial D$
such that $w$ belongs to the arc in $\Gamma $ between $w^{\prime }$ and $%
w^{\prime \prime }$ and apart from its end points, this arc is contained in $%
D$. And this is a contradiction with the Poincar\'e recurrence theorem, see
the end of the proof of lemma \ref{nodisctrans}. $\Box $

\vskip 0.2truecm

Back to the proof of theorem \ref{casogeral}, as $\widetilde{f}^m(\overline{%
W^u(\widetilde{Q}_0)})=\overline{W^u(\widetilde{Q}_0)}$ and this set is
closed and connected, any connected component $\widetilde{M}$ of its
complement is an open disk such that $\widetilde{M}\cap (\widetilde{M}%
+(a,b))=\emptyset $ for all integers $(a,b)\neq (0,0).$ To see this, suppose
by contradiction that for some $(a^{\prime },b^{\prime })\neq (0,0),$ $
\widetilde{M}$ intersects $\widetilde{M}+(a^{\prime },b^{\prime }).$ As all
integer translates of $\widetilde{M}$ are in the complement of $\overline{%
W^u(\widetilde{Q}_0)},$ this contradicts the fact that $\overline{W^u( 
\widetilde{Q}_0)}$ is invariant under integer translations and connected. So 
$M=p(\widetilde{M})$ is a open disk in the torus and there exists an integer 
$k>0$ such that 
$$
f^{k.m}(M)=M. 
$$
This follows from the fact that $M$ is a connected component of the
complement of $\overline{W^u(Q_0)},$ which is invariant under $f^m.$ By
lemma \ref{discgeral}, $diam(\widetilde{M})$ is uniformly bounded from
above. An analogous argument applied to $\overline{W^s(\widetilde{Q}_0)}$
implies that any connected component of the complement of $\overline{W^s( 
\widetilde{Q}_0)}$ is also a connected component of the lift of an open disk
of the torus, with diameter uniformly bounded from above. 

In order to prove that $\overline{W^u(\widetilde{Q}_0)}=\overline{W^s( 
\widetilde{Q}_0)},$ from expression (\ref{ulttimma}) it is enough to show
that $\overline{W^u(Q_0)}=\overline{W^s(Q_0)}.$ If $M$ is a connected
component of the complement of $\overline{W^s(Q_0)},$ from what we did
above, there exists an integer $k>0$ such that $f^{k.m}(M)=M$ and $M$ is an
open disk in the torus, with diameter bounded from above.

\begin{description}
\item[Claim 1]  : $W^u(Q_0)$ does not intersect $M.$
\end{description}

{\it Proof:}

If $W^u(Q_0)$ intersects $M,$ as any connected component of $p^{-1}(M)$ has
the same diameter, which is bounded from above, and for all pair of integers 
$(i,j),$ $W^u(\widetilde{Q}_0)\pitchfork W^s(\widetilde{Q}_0+(i,j)),$ we get
that some branch $\Gamma $ of $W^u(Q_0)$ intersects $M$ at some point $w$
and there are points $w^{\prime },w^{\prime \prime }\in \Gamma \cap \partial
M$ such that $w$ belongs to the arc in $\Gamma $ between $w^{\prime }$ and $%
w^{\prime \prime }$ and apart from its end points, this arc is contained in $%
M.$ As in the end of the proof of lemma \ref{discgeral}, this is a
contradiction with the Poincar\'e recurrence theorem. $\Box $

\vskip 0.2truecm

So $\overline{W^u(Q_0)}\subset \overline{W^s(Q_0)}$ and an analogous
argument gives $\overline{W^s(Q_0)}\subset \overline{W^u(Q_0)},$ which
finishes this part of the argument.

Finally, let us prove that 
\begin{equation}
\label{invariance}\widetilde{f}(\overline{W^u(\widetilde{Q}_0)})=\overline{%
W^u(\widetilde{Q}_0)}. 
\end{equation}
The equalities in expression (\ref{ulttimma}) imply that expression (\ref
{invariance}) is true if and only if, 
\begin{equation}
\label{invariancetorus}f(\overline{W^u(Q_0)})=\overline{W^u(Q_0)}. 
\end{equation}
If expression (\ref{invariancetorus}) does not hold, then there exists a
connected component $M$ of $\left( \overline{W^u(Q_0)}\right) ^c$ such that
either:

\begin{enumerate}
\item  $W^u(f(Q_0))$ intersects $M;$

\item  $W^u(Q_0)$ intersects $f(M);$
\end{enumerate}

As $M$ is a bounded $f$-periodic open disk, claim 1 implies that both cases
above can not happen. $\Box $

\vskip 0.2truecm

\subsection{Proof of lemma 2}

{\bf Statement}{\it :}\ {\it Suppose $f$ belongs to $Diff_0^{1+\epsilon }(%
{\rm T^2})$ and $int(\rho (\widetilde{f}))$ is not empty or $f$ belongs to $%
Diff_k^{1+\epsilon }({\rm T^2})$ and $int(\rho _V(\widetilde{f}))$ is not
empty. If $f$ is transitive, then $f$ can not have a periodic open disk. In
the general case, there exists $M=M(f)>0$ such that if $D\subset {\rm T^2}$
is a $f$-periodic open disk, then for any connected component $\widetilde{D}$
of $p^{-1}(D),$ $diam(\widetilde{D})<M.$}

\vskip 0.2truecm

This is contained in lemmas \ref{nodisctrans} and \ref{discgeral}. $\Box $

\vskip 0.2truecm

\subsection{Proof of corollary 7}

{\bf Statement}{\it :}\ {\it Suppose $f$ belongs to $Diff_0^{1+\epsilon }(%
{\rm T^2})$ and $(0,0)\in int(\rho (\widetilde{f}))$ or $f$ belongs to $%
Diff_k^{1+\epsilon }({\rm T^2})$ and $0\in int(\rho _V(\widetilde{f})).$ If $
\widetilde{D}\subset {\rm I}\negthinspace {\rm R^2}$ is a wandering open
disk, then for all integers $n,$ $diam(\widetilde{f}^n(\widetilde{D}))$ is
uniformly bounded from above and there exists a $f$-periodic open disk $%
D_{+}\subset {\rm T^2}$ such that $D_{+}\supset D=p(\widetilde{D}).$}

\vskip 0.2truecm

From theorem \ref{casogeral}, $\widetilde{f}$ has a hyperbolic $m$-periodic
point $\widetilde{Q}$ such that for all pair of integers $(a,b),$ $\overline{%
W^u(\widetilde{Q})}=\overline{W^s(\widetilde{Q})}=\overline{W^u(\widetilde{Q}%
)}+(a,b)$ and all connected components of $\left( \overline{W^u(\widetilde{Q}%
)}\right) ^c$ are open disks, with uniformly bounded diameter.

So, if $diam(\widetilde{f}^n(\widetilde{D}))$ is not uniformly bounded, then
for some integer $n_0,$ $\widetilde{f}^{n_0}(\widetilde{D})$ intersects both 
$W^s(\widetilde{Q})$ and $W^u(\widetilde{Q}).$ But this means that there
exists a subset of $\widetilde{f}^{n_0-n.m}(\widetilde{D})$ that gets closer
and closer (in the Hausdorff topology) as $n\rightarrow \infty ,$ to a piece
of $W^s(\widetilde{Q})$ that contains $\widetilde{Q}$ and in a similar way,
there exists a subset of $\widetilde{f}^{n_0+n.m}(\widetilde{D})$ that gets
closer and closer (also in the Hausdorff topology) as $n\rightarrow \infty ,$
to a piece of $W^u(\widetilde{Q})$ that also contains $\widetilde{Q}.$ Thus
if $n>0$ is sufficiently large, $\widetilde{f}^{n_0-n.m}(\widetilde{D})$
intersects $\widetilde{f}^{n_0+n.m}(\widetilde{D})$, which implies that 
$$
\widetilde{f}^{2n.m}(\widetilde{D})\text{ intersects }\widetilde{D}, 
$$
a contradiction with the assumptions on $\widetilde{D}.$ In a similar way,
if $\widetilde{D}$ intersects $W^u(\widetilde{Q}),$ then theorem \ref
{casogeral} implies that $\widetilde{D}$ intersects $W^s(\widetilde{Q})$ and
so it can not be wandering. Thus $\widetilde{D}\subset \left( \overline{W^u( 
\widetilde{Q})}\right) ^c$ and the corollary follows by noticing that $
\widetilde{D}_{+}$ is the connected component of $\left( \overline{W^u( 
\widetilde{Q})}\right) ^c$ that contains $\widetilde{D}$ and $D_{+}=p( 
\widetilde{D}_{+}).$ $\Box $

\vskip 0.2truecm

\subsection{Proof of corollary 8}

{\bf Statement}{\it :}\ {\it If $f$ belongs to $Diff_0^{1+\epsilon }({\rm T^2%
})$ and $\{(\frac pq,\frac rq),(0,0)\}\in int(\rho (\widetilde{f}))$ then $
\widetilde{f}^q(\bullet )-(p,r)$ has a hyperbolic periodic point $\widetilde{%
P}$ such that for any pair of integers $(a,b),$ $W^u(\widetilde{P})%
\pitchfork W^s(\widetilde{P}+(a,b))$ and $R.I.(f)=\overline{W^u(p(\widetilde{%
P}))}.$ In case $f$ belongs to $Diff_k^{1+\epsilon }({\rm T^2})$ and $%
\{\frac pq,0\}\in int(\rho _V(\widetilde{f})),$ then, for all integers $s,$ $
\widetilde{f}^q(\bullet )-(s,p)$ has a hyperbolic periodic point $\widetilde{%
P}$ such that for any pair of integers $(a,b),$ $W^u(\widetilde{P})%
\pitchfork W^s(\widetilde{P}+(a,b))$ and $R.I.(f)=\overline{W^u(p(\widetilde{%
P}))}.$ Moreover, if $\widetilde{Q}$ is the $\widetilde{f}$-periodic point
that appears in the statement of theorem \ref{casogeral}, then $W^u(p( 
\widetilde{P}))\pitchfork W^s(p(\widetilde{Q}))$ and $W^u(p(\widetilde{Q}))%
\pitchfork W^s(p(\widetilde{P})).$}

\vskip 0.2truecm

The complement of $R.I.(f)$ is a union of $f$-periodic disks (with diameter
bounded from above when lifted to the plane) and by theorem \ref{casogeral},
there exists a point $\widetilde{P}\in {\rm I}\negthinspace {\rm R^2}$ such
that $p(\widetilde{P})$ is a hyperbolic $f$-periodic point with adequate
rotation vector or number and for any pair of integers $(a,b),$ $W^u( 
\widetilde{P})\pitchfork W^s(\widetilde{P}+(a,b)).$ Moreover, the complement
of $\overline{W^u(p(\widetilde{P}))}$ is also a union of $f$-periodic disks,
with diameter uniformly bounded from above when lifted to the plane. So if $%
\left( R.I.(f)\right) ^c\neq \left( \overline{W^u(p(\widetilde{P}))}\right)
^c,$ then there exists a $f$-periodic disk $D$ which intersects the unstable
manifold of a hyperbolic $f$-periodic point $R$ ($R$ can be equal to {\it $%
p( \widetilde{Q})$} or {\it $p(\widetilde{P})$}) such that the union of its
stable and unstable manifolds contains closed curves in the torus of all
homotopy types. This implies that some branch $\Gamma $ of $W^s(R)$ or of $%
W^u(R)$ intersects $D$ at some point $w$ and there are points $w^{\prime
},w^{\prime \prime }\in \Gamma \cap \partial D$ such that $w$ belongs to the
arc in $\Gamma $ between $w^{\prime }$ and $w^{\prime \prime }$ and apart
from its end points, this arc is contained in $D.$ As in the end of the
proof of lemma \ref{discgeral}, this is a contradiction with the Poincar\'e
recurrence theorem. So $\left( R.I.(f)\right) ^c=\left( \overline{W^u(p( 
\widetilde{P}))}\right) ^c,$ which implies that $R.I.(f)=\overline{W^u(p( 
\widetilde{P}))}.$ Now we prove the last part of the corollary.

Let us choose a closed curve $\gamma \subset ${\it ${\rm T^2,}$} homotopic
to $(1,0),$ such that $\gamma $ is the projection under $p$ of a simple
connected arc $\widetilde{\gamma }\subset ${\it ${\rm I}\negthinspace {\rm %
R^2}$} made of two connected pieces:\ one starts at $\widetilde{Q},$
contained in $W^u(\widetilde{Q})$ until it reaches $W^s(\widetilde{Q})+(1,0)$
and the other piece is contained in $W^s(\widetilde{Q})+(1,0)$ until it
reaches $\widetilde{Q}+(1,0).$ As {\it $p(\widetilde{P})\neq p(\widetilde{Q}%
),$} we can choose a closed curve $\alpha \subset W^s(p(\widetilde{P}))\cup
p(\widetilde{P})\cup W^u(p(\widetilde{P}))$ whose homotopy class is $(0,1)$,
such that $\alpha \cap W^s(p(\widetilde{P}))$ is a connected arc,
sufficiently close to $p(\widetilde{P}),$ in a way that $\gamma \cap \left(
\alpha \cap W^s(p(\widetilde{P}))\right) =\emptyset.$ This is possible
because $p(\widetilde{P})\cap \gamma=\emptyset$.

As $\alpha $ and $\gamma $ must intersect, the previous choice made on $%
\alpha $ implies that {\it $W^u(p(\widetilde{P}))\pitchfork W^s(p(\widetilde{%
Q})).$} Interchanging $\widetilde{Q}$ and $\widetilde{P}$ in the above
argument we also obtain that {\it $W^u(p(\widetilde{Q}))\pitchfork W^s(p( 
\widetilde{P}))$} and the proof is over. $\Box $

\vskip 0.2truecm

\subsection{Proof of theorem 9}

{\bf Statement}{\it :} {\it Suppose $f$ belongs to $Diff_0^{1+\epsilon }(%
{\rm T^2})$ and $(0,0)\in int(\rho (\widetilde{f}))$ or $f$ belongs to $%
Diff_k^{1+\epsilon }({\rm T^2})$ and $0\in int(\rho _V(\widetilde{f})).$
Then, for all $\theta \in [0,2\pi ],$ $\omega (B_\theta ^\infty )=\emptyset $
and $\omega (B_S^\infty )=\omega (B_N^\infty )=\emptyset .$}

\vskip 0.2truecm

The case when $\theta $ is a rational multiple of $2\pi ,$ $f\in
Diff_0^{1+\epsilon }({\rm T^2})$ and $\left( 0,0\right) \in int(\rho ( 
\widetilde{f}))$ is analogous to the case when $f\in Diff_k^{1+\epsilon }(%
{\rm T^2})$ and $0\in int(\rho _V(\widetilde{f})).$ In this case, if $f\in
Diff_0^{1+\epsilon }({\rm T^2}),$ by conjugating $f$ with a suitable map, we
can assume that $\theta =0.$ So, suppose by contradiction that $\omega
(B_0^\infty )\neq \emptyset .$ Let $\Gamma $ be a connected component of $%
\omega (B_0^\infty ).$ Then $\Gamma \subset [0,\infty [\times {\rm I}%
\negthinspace 
{\rm R}$ and $\Gamma $ is closed and unbounded. %
%
Now, we remember theorem \ref{casogeral} and let $\widetilde{Q}$ be a
hyperbolic $m$-periodic point for $\widetilde{f}$ such that for any pair of
integers $(a,b),$ $W^u(\widetilde{Q})\pitchfork W^s(\widetilde{Q}+(a,b))$
and $\overline{W^u(\widetilde{Q})}=\overline{W^s(\widetilde{Q})}.$ So we can
pick two simple arcs, $\alpha _{(1,0)}$ and $\alpha _{(0,1)}$ such that $%
\alpha _{(1,0)}\subset W^u(\widetilde{Q})\cup W^s(\widetilde{Q}+(1,0))$ and $%
\alpha _{(1,0)}$ connects $\widetilde{Q}$ to $\widetilde{Q}+(1,0)$ and
similarly, $\alpha _{(0,1)}\subset W^u(\widetilde{Q})\cup W^s(\widetilde{Q}%
+(0,1))$ and $\alpha _{(0,1)}$ connects $\widetilde{Q}$ to $\widetilde{Q}%
+(0,1).$ Let us define two path connected sets $\theta _{(1,0)},\theta
_{(0,1)}\subset {\rm I}\negthinspace {\rm R^2}$ in the following way: $%
\theta _{(1,0)}=\stackunder{i\in integers}{\cup }\alpha _{(1,0)}+i(1,0)$ and 
$\theta _{(0,1)}=\stackunder{i\in integers}{\cup }\alpha _{(0,1)}+i(0,1).$
Then, the next consequences hold:

\begin{enumerate}
\item  $\theta _{(1,0)}=\theta _{(1,0)}+(1,0)$ and $\theta _{(0,1)}=\theta
_{(0,1)}+(0,1);$

\item  $\theta _{(1,0)}$ contains $\widetilde{Q}+i(1,0)$ and $\theta
_{(0,1)} $ contains $\widetilde{Q}+i(0,1),$ for all integers $i;$
\end{enumerate}

It is easy to see that if $\theta _{(0,1)}+i(1,0)$ does not intersect $%
\Gamma $ for all integers $i,$ then as $\Gamma $ is connected and unbounded, 
$\theta _{(1,0)}+j(0,1)$ must intersect $\Gamma $ for some integer $j.$ From
this we get that either $W^u(\widetilde{Q}+(c,d))$ or $W^s(\widetilde{Q}%
+(c,d)),$ for some pair of integers $(c,d),$ has a topologically transverse
intersection with $\Gamma .$ Suppose it is $W^u(\widetilde{Q}+(c,d)).$ By
theorem \ref{casogeral} we can suppose that $(c,d)=(0,0).$ This implies that 
$\widetilde{f}^{-n.m}(\Gamma )$ contains a compact connected set that
converges in the Hausdorff topology to $W^s(\widetilde{Q})$ as $n\rightarrow
\infty ,$ a contradiction because $\widetilde{f}(\omega (B_0^\infty
))=\omega (B_0^\infty )\subset [0,\infty [\times {\rm I}\negthinspace 
{\rm R}$ and $W^s(\widetilde{Q})$ has points outside $[0,\infty [\times {\rm %
I}\negthinspace 
{\rm R.}$

In case $\frac \theta {2\pi }$ is irrational, we proceed as follows: The
theory of continued fractions (see for instance \cite{kinchin}) implies that
using integer translates of the arcs $\alpha _{(1,0)}$ and $\alpha _{(0,1)}$
defined above, we can obtain a connected closed set $\gamma $ in the plane
such that $\gamma $ is contained between two straight lines of slope $\theta 
$ and $\gamma $ intersects all straight lines of slope $\theta +\pi /2.$ We
can do the same procedure for $\theta +\pi /2,$ that is, also using integer
translates of $\alpha _{(1,0)}$ and $\alpha _{(0,1)}$ we can obtain a
connected closed set $\gamma ^{*}$ in the plane such that $\gamma ^{*}$ is
contained between two straight lines of slope $\theta +\pi /2$ and $\gamma
^{*}$ intersects all straight lines of slope $\theta .$ So if $\omega
(B_\theta ^\infty )\neq \emptyset ,$ then some integer translate of $\gamma
\cup \gamma ^{*}$ intersects $\omega (B_\theta ^\infty ).$ From this and
theorem \ref{casogeral}, we get that either $W^u(\widetilde{Q})$ or $W^s( 
\widetilde{Q})$ has a topologically transverse intersection with some
connected component $\Gamma $ of $\omega (B_\theta ^\infty ).$ Suppose it is 
$W^u(\widetilde{Q}).$ As in the rational case, this implies that $\widetilde{%
f}^{-n.m}(\Gamma )$ contains a compact connected set that converges in the
Hausdorff topology to $W^s(\widetilde{Q})$ as $n\rightarrow \infty ,$ a
contradiction because $\widetilde{f}(\omega (B_\theta ^\infty ))=\omega
(B_\theta ^\infty )\subset \{\widetilde{z}\in {\rm I}\negthinspace 
{\rm R^2:}\left\langle \widetilde{z},(\cos (\theta ),\sin (\theta
))\right\rangle \geq 0\}$ and $W^s(\widetilde{Q})$ has points outside this
set. $\Box $

\vskip 0.2truecm

\subsection{Proof of corollary 10}

{\bf Statement}{\it :}\ {\it Suppose $f$ belongs to $Diff_0^{1+\epsilon }(%
{\rm T^2})$ and $(0,0)\in int(\rho (\widetilde{f}))$ or $f$ belongs to $%
Diff_k^{1+\epsilon }({\rm T^2})$ and $0\in int(\rho _V(\widetilde{f})).$ If $%
f$ is transitive, then $\overline{p(B_S^\infty )}=\overline{p(B_N^\infty )}= 
\overline{p(B_\theta ^\infty )}={\rm T^2,}$ for all $\theta \in [0,2\pi ].$
And for a general $f,$ any of the following sets $\left( \overline{%
p(B_S^\infty )}\right) ^c,$ $\left( \overline{p(B_N^\infty )}\right) ^c$ and 
$\left( \overline{p(B_\theta ^\infty )}\right) ^c$ is the union of $f$%
-periodic open disks, with diameters uniformly bounded from above when
lifted to the plane.}

\vskip 0.2truecm

Without loss of generality, suppose $f$ belongs to $Diff_0^{1+\epsilon }(%
{\rm T^2})$ and $(0,0)\in int(\rho (\widetilde{f})).$ Given $\theta \in
[0,2\pi ],$ let us look at a connected component $U$ of $\left( \overline{%
p(B_\theta ^\infty )}\right) ^c.$ Note that $U$ is periodic, because $f( 
\overline{p(B_\theta ^\infty )})=\overline{p(B_\theta ^\infty )}.$ %
%
%
%
So there are two possibilities:

\begin{enumerate}
\item  $U$ contains a homotopically non-trivial simple closed curve $\gamma $
in the torus;

\item  $U$ does not contain such a curve;
\end{enumerate}

In the first case, without loss of generality, suppose $\gamma $ is a
vertical curve, that is $\gamma $ is homotopic to $(0,1).$ As $(0,0)\in
int(\rho (\widetilde{f})),$ if $\widetilde{\gamma }$ is a connected
component of $p^{-1}(\gamma ),$ then for all integers $n,$ $\widetilde{f}^n( 
\widetilde{\gamma })$ intersects $\widetilde{\gamma },$ otherwise $(0,0)$
would not be an interior point of $\rho (\widetilde{f}).$ And given any
integer $i\neq 0,$ we get that for some integer $n(i)>0,$ $\widetilde{f}%
^{n(i)}(\widetilde{\gamma })$ intersects $\widetilde{\gamma }+(i,0).$ So
some vertical translate of $B_\theta ^\infty $ intersects $\widetilde{f}%
^{n(2)}(p^{-1}(\gamma )),$ which means that $p(B_\theta ^\infty )$
intersects $f^{n(2)}(\gamma ),$ a contradiction. Thus case 1 does not happen.

In the second case let us prove that $U$ is an open disk. If it is not, then
there exists a simple closed curve $\alpha \subset U,$ which is contractible
as a curve in the torus, such that in the disk bounded by $\alpha $ there
are points of $\overline{p(B_\theta ^\infty )}.$ But this is impossible,
because each connected component $\Gamma $ of $B_\theta ^\infty $ is an
unbounded closed subset of the plane. So $U$ is an open disk.

If $f$ is transitive, lemma \ref{nodisctrans} tells us that case 2 also does
not happen, which means that $\overline{p(B_\theta ^\infty )}={\rm T^2.}$ In
the general case, $\left( \overline{p(B_\theta ^\infty )}\right) ^c$ is
equal the union of $f$-periodic open disks, which by lemma \ref{discgeral},
all have diameters uniformly bounded from above when lifted to the plane. $%
\Box $

\vskip 0.2truecm

{\it Acknowledgements: }I would like to thank Andr\'e de Carvalho for
telling me about his results in \cite{andre} and for his help during the
preparation of this paper. I am also very grateful to the referee for all
his comments and suggestions.

\centerline{\bf Figure captions.}

\begin{itemize}
\item[Figure 1. ]  Diagram showing a closed connected set $K$ which
satisfies $K\pitchfork W^u(P).$

\item[Figure 2.]  Diagram showing the simple arc $\eta .$

\item[Figure 3.]  Diagram showing the simple arc $\eta $ and how the set $%
\theta $ is obtained.

\item[Figure 4.]  Diagram showing the points $\widetilde{z},\widetilde{Q}%
_{V1},\widetilde{Q}_{V2}$ and the rectangle $\widetilde{Ret}(\widetilde{z},
\widetilde{Q}_{V1},\widetilde{Q}_{V2}).$
\end{itemize}

\begin{center}
\mbox{\includegraphics[width=13cm]{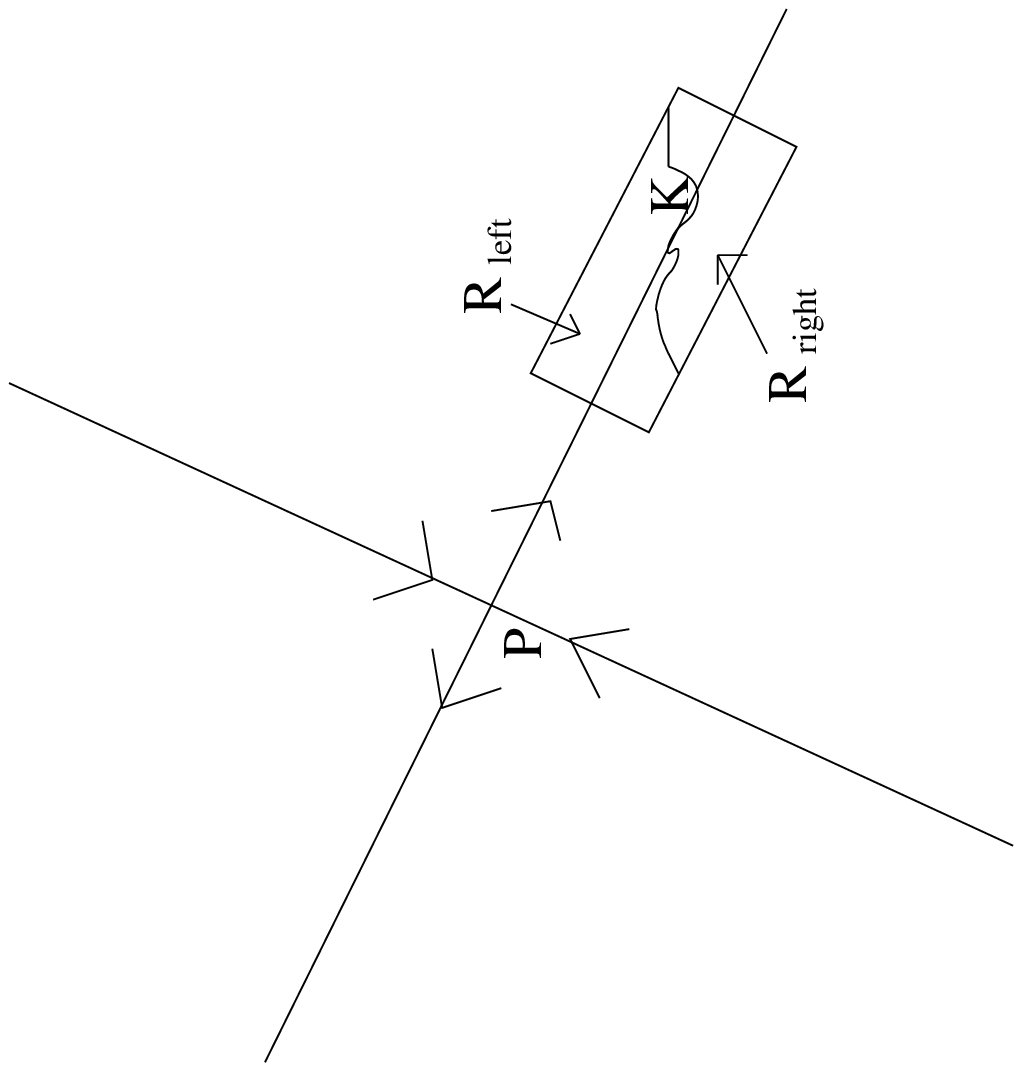}} 
\mbox{\includegraphics[width=13cm]{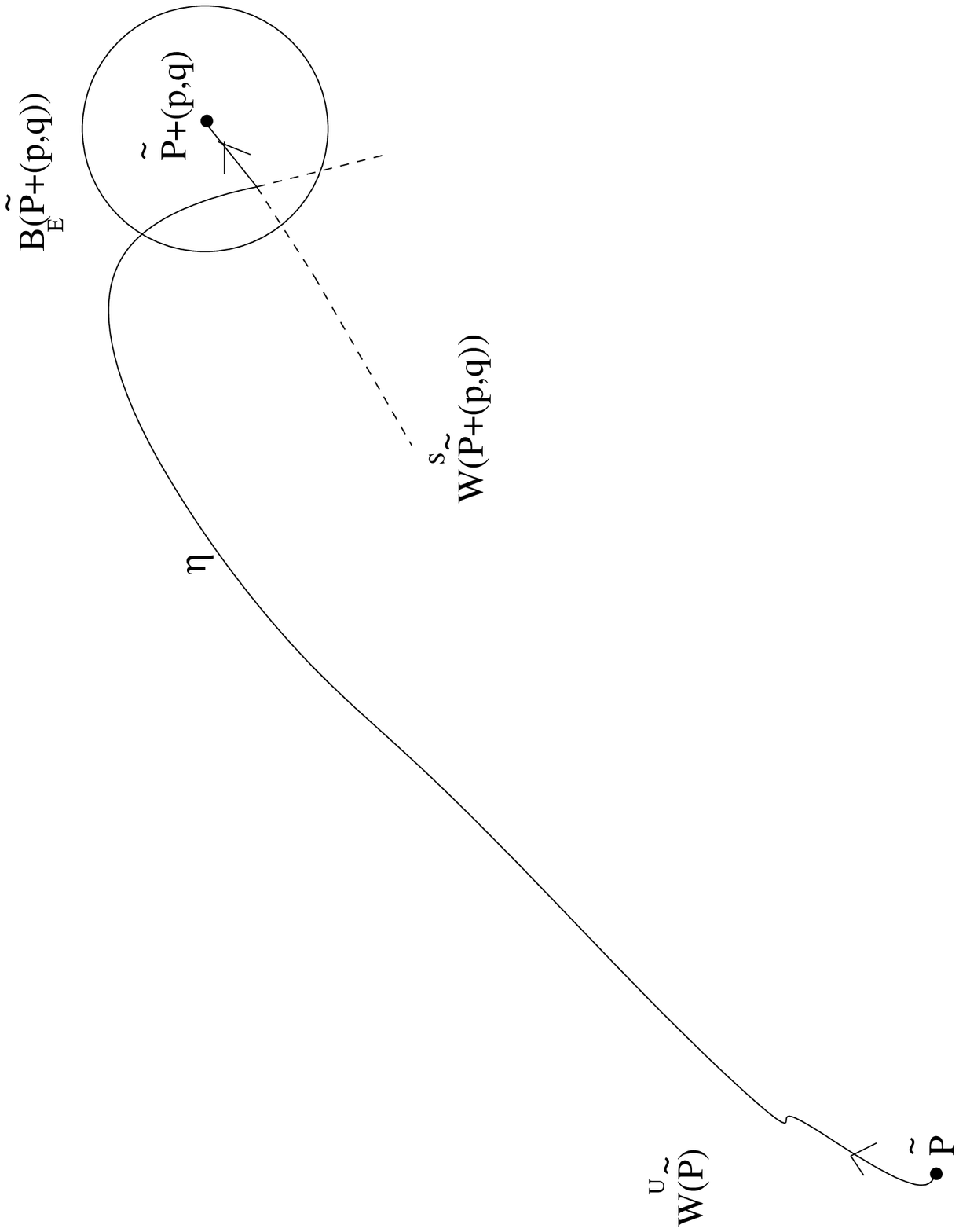}} 
\mbox{\includegraphics[width=13cm]{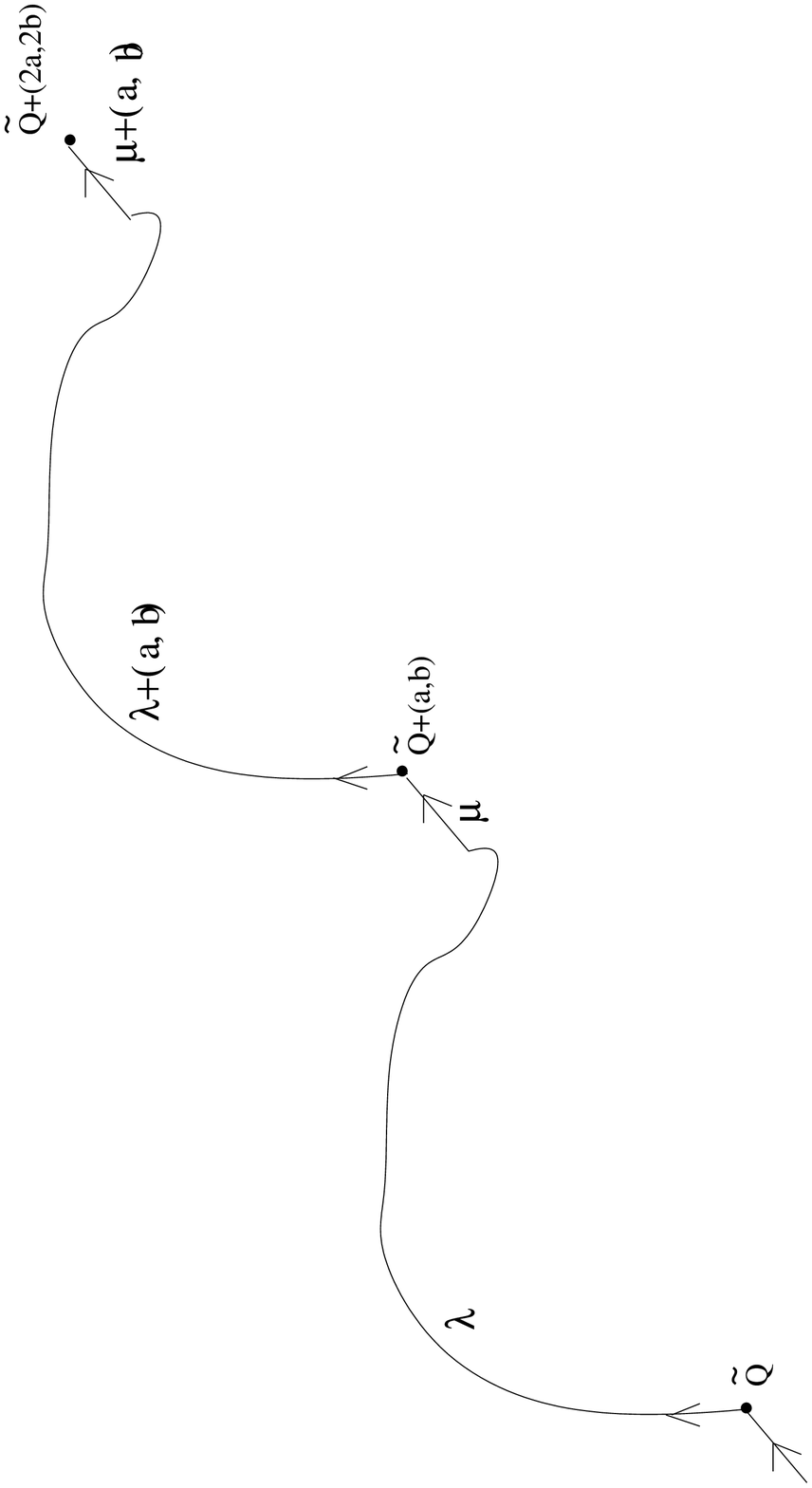}} 
\mbox{\includegraphics[width=13cm]{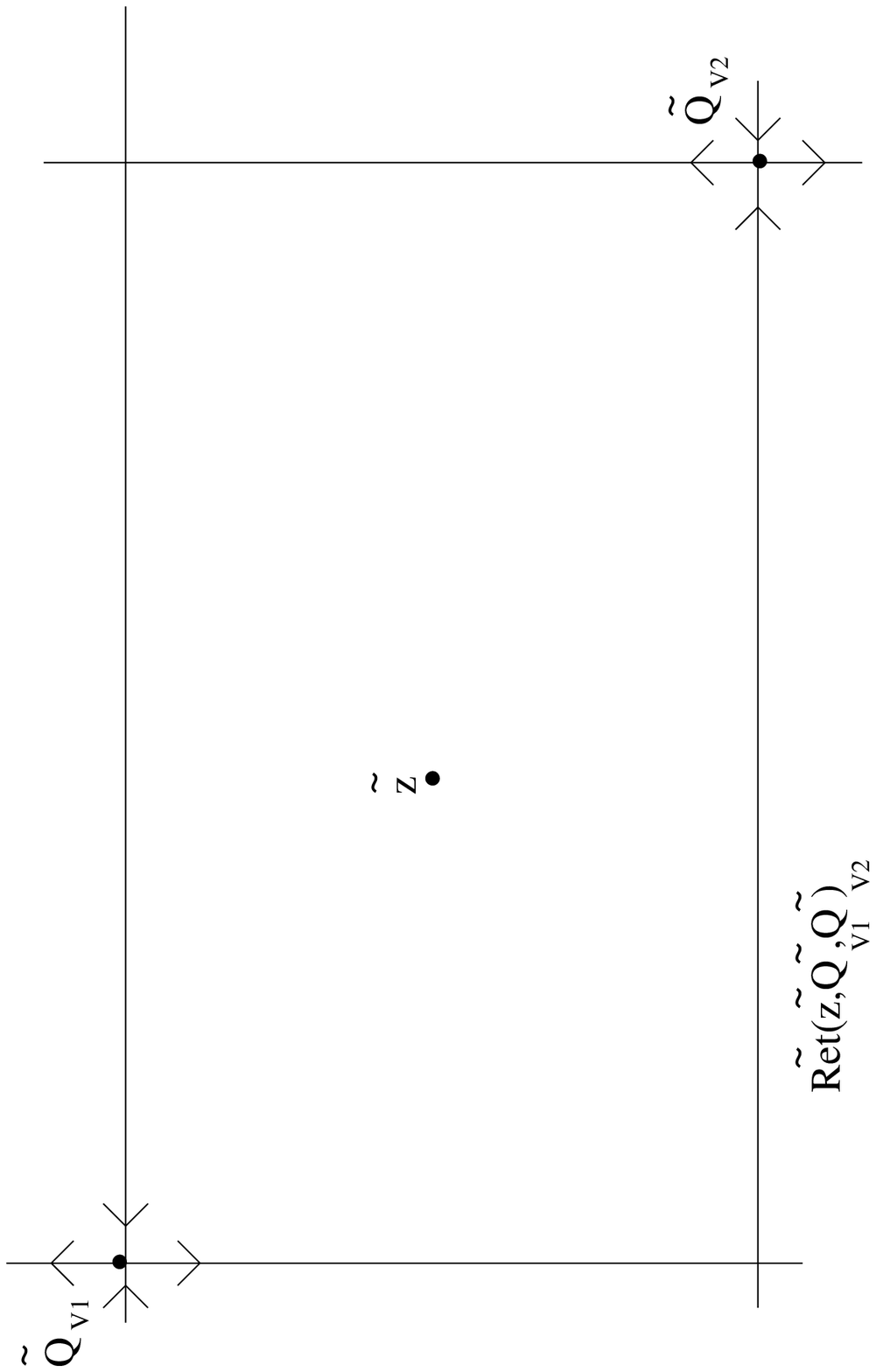}}
\end{center}

\end{document}